\documentclass[10pt,leqno]{amsart}
\usepackage{graphicx}

\usepackage{indentfirst,csquotes}

\usepackage{cite}
\usepackage{amsmath,amssymb,amsfonts}
\usepackage{algorithmic}
\usepackage{textcomp}
\usepackage{xcolor}
\usepackage{braket}
\usepackage{comment}
\usepackage{mathdots}
\usepackage{paralist,hyperref,etoolbox}
\usepackage{lipsum}
\usepackage{mathrsfs}

\usepackage{todonotes}

\usepackage{amssymb,amsthm,amsmath}

\newtheorem{theorem}{Theorem}

\newtheorem{lemma}{Lemma}

\newtheorem{corollary}{Corollary}

\newtheorem{remark}{Remark}

\hypersetup{
    colorlinks=true,
    linkcolor=black,
    filecolor=black,
    urlcolor=black
}

\begin{document}
\title[LHAM to Solve Nonlinear PDEs]{Lindbladian Homotopy Analysis Method to Solve Nonlinear Partial Differential Equations} 
\author[Choi, et al.]{
Eunsik Choi$^{1}$, Jungin E. Kim$^{1}$, Xueling Lu$^{2}$, Yan Wang$^{1}$\\[1ex]
$^{1}$Georgia Institute of Technology, Atlanta, GA, USA\\
$^{2}$Qscico, Atlanta, GA, USA\\[1ex]
}

\date{\today}

\maketitle

\let\thefootnote\relax

\begin{abstract}
Quantum scientific computing is to solve engineering and science problems such as simulation and optimization on quantum computers. Solving ordinary and partial differential equations (PDEs) is essential in simulations. However, existing quantum approaches to solve nonlinear PDEs suffer from the issues of curse of dimensionality and convergence during the linearization process.
In this paper, a Lindbladian homotopy analysis method (LHAM) is proposed as a quantum differential equation solver to simulate nonlinear and nonunitary dynamics.
The original nonlinear problem is first converted to a recursive sequence of linear PDEs with the homotopy analysis method and reformulated as a higher-dimensional lower block triangular linear homogeneous autonomous system. 
The solution is then embedded in the density matrix and obtained through the Lindblad dynamics simulation. 
Compared to other methods such as Carleman linearization and the Koopman-von Neumann approach where the dimension of Hilbert space increases polynomially with the inverse of truncation error, the Hilbert space dimension in LHAM increases only logarithmically.
LHAM is demonstrated with nonlinear PDEs including Burgers' equation and reduced magnetohydrodynamics equations. 
\end{abstract} 

\bigskip

\noindent 
\section{\label{sec:intro}Introduction}

Quantum scientific computing is to solve engineering and science problems such as simulation and optimization on quantum computers. 
Particularly, simulations of complex phenomena such as fluid dynamics, phase transitions, shallow water waves, and plasma dynamics require solving nonlinear ordinary differential equations (ODEs) and partial differential equations (PDEs). 
It is challenging to solve nonlinear ODEs and PDEs on quantum computers because quantum computing relies on linear unitary operations.

Different quantum algorihtms have been developed to solve nonlinear ODEs and PDEs, such as computing multiple copies of states \cite{leyton2008quantum}, nonlinear Schr\"odinger linearization \cite{lloyd2020quantum}, 
the Koopman-von Neumann approach \cite{joseph2020koopman, dodin2021applications, giannakis2022embedding, jin2022quantum, jin2023time}, Carleman linearization \cite{engel2021linear, liu2021efficient, akiba2023carleman, itani2024quantum, wu2025quantum, brustle2025quantum, gonzalez2025quantum}, 
and variational quantum algorithms \cite{lubasch2020variational, kyriienko2021solving, jaksch2023variational, su2025quantum}.
Most of the above linearization methods rely on increasing the dimension of the solution space associated with the original nonlinear problem.
As a result, they suffer from the issues of curse of dimensionality and convergence.

Recently, the homotopy analysis method \cite{liao2004homotopy} was adopted to solve nonlinear ODEs and PDEs on quantum computers \cite{xue2021quantum, xue2025quantum, bharadwaj2025quantum, liao2026ham}.
The homotopy analysis method alleviates both the dimensionality and convergence issues in linearization.
With the introduction of an embedding parameter, a homotopy function is defined so that it can be expanded as a homotopy-Maclaurin series with respect to the initial guess of nonlinear PDE solution. The homotopy-Maclaurin series is used to obtain a recursive system of linear PDEs known as deformation equations. 
The nonhomogeneous term of each linear PDE only depends on the solutions of the lower-order deformation equations.
Thus, the problem of solving a nonlinear PDE is converted to solving multiple nonhomogeneous linear PDEs without increasing the dimension of the solution space.

After linearization, linear ODEs and PDEs can be solved with different methods. One method is quantum linear equation solver \cite{berry2014high}, where the original differential equation is discretized over time and then solved with a quantum linear system algorithm \cite{harrow2009quantum}.
The complexity of quantum linear equation solver depends on condition number, which is particularly large for engineering problems.
Another method is variational quantum simulation \cite{endo2020variational, sul2024quantum, sul2025generic}, where the quantum circuit parameters are evolved over time by solving McLachlan's variational principle. 
However, its accuracy is limited when the exact trajectory of time evolution is not projected onto the tangent space of the variational manifold corresponding to the parametric ansatz.
The third method is variational quantum linear solver \cite{bravo2023variational} that optimize parameterized circuits to approximate the evolving state or the solution of a linear system. It requires heuristic optimization of parameters.
The fourth method is Schr\"odingerization \cite{jin2023quantum, jin2024quantum} or linear combination of Hamiltonian simulation \cite{an2023linear}, where the general dynamics is converted to unitary dynamics by introducing an auxiliary variable. 
However, the quantum circuits are implemented with ancilla qubits which caused the Hilbert space dimension to increase. 
The number of ancilla qubits increases logarithmically as the resolution of the discretized auxiliary variable increases. 
As the number of ancilla qubits increases, the quantum circuit depth increases exponentially.

In this paper, a Lindbladian homotopy analysis method (LHAM) is proposed to simulate nonlinear and nonunitary dynamics. 
In LHAM, a nonlinear PDE is linearized with the homotopy analysis method. 
The homotopy analysis method avoids the exponential increase of the state space dimension which occurs in Carleman linearization and the Koopman-von Neumann approach.
The second contribution of this paper is the introduction of new auxiliary variables to represent the dynamics of the nonhomogeneous forcing term in deformation equations. The auxiliary variables couple solutions of different homotopy orders through the off-diagonal feed-forward coupling operator so that the homotopy-series solution can be obtained by solving one linear homogeneous autonomous system.
The set of the auxiliary variables is closed under the differentiation operation by introducing time monomial factor when the homogeneous and nonhomogeneous terms are resonant with the same decay rate. 
The third contribution is the incorporation of mid-circuit-measurement-and-reset strategy and high-order Stinespring dilation in simulating Lindblad-based nonunitary dynamics for solving  differential equations.

The nonunitary dynamics simulation method of LHAM is inspired by the recently developed Lindbladian embedding \cite{shang2025designing}, where the solution of the linear ODE is encoded in the off-diagonal entries of the density matrix.
The simulation of the nonunitary dynamics is achieved by embedding system evolutions in quantum channels and solving the Lindblad master equation. 
An ancilla qubit is used to encode Hermitian and anti-Hermitian components of the system operator.
Our extension from Ref. \cite{shang2025designing} is an efficient implementation based on the Stinespring dilation, where as few as one environment ancilla qubit is needed with the mid-circuit-measure-and-reset strategy.
LHAM only requires as few as two ancilla qubits to simulate nonunitary dynamics. 
We further generalize and extend the strategy by introducing higher-order Stinespring dilations for more accurate simulations.

The proposed LHAM exhibits improved scalability in comparison with the Carleman linearization and the Koopman-von Neumann approach.
In LHAM, the dimension of the Hilbert space scales in the order of $\mathscr{O}(D\tilde{m})$, where $D$ is the dimension of the solution space after discretization and $\tilde{m}$ is the homotopy truncation order. 
$\tilde{m}$ scales in the order of $\mathscr O(\mathrm{log}(1/\varepsilon))$, where $\varepsilon$ is the truncation error \cite{hetmaniok2014usage}.
Thus, the dimension of the Hilbert space in LHAM increases with $\mathscr{O}(D\mathrm{log}(1/\varepsilon))$.
In contrast, the dimension in Carleman linearization scales as $\mathscr{O}(D^N)$ \cite{engel2021linear},
where $N$ is the maximum polynomial degree
and scales as $N \sim \mathscr O(\mathrm{log}(1/\varepsilon))$ \cite{wu2025quantum}. 
Thus, the dimension in Carleman linearization increases polynomially with $1/\varepsilon$ as $\mathscr{O}((1/\varepsilon)^{\mathrm{log}D})$.
In the Koopman-von Neumann approach, $D$ denotes the dimension of the nonlinear ODE system obtained after spatial discretization of the original nonlinear PDE. The Koopman-von Neumann wavefunction is defined over the resulting $D$-dimensional continuous state space. If each state-space coordinate is discretized using $N$ grid points, the discretized Koopman-von Neumann Hilbert space has dimension $\mathscr O(N^D)$ \cite{joseph2020koopman}, where $N \sim \mathscr{O}(1/\varepsilon)$ \cite{jin2023time}. 
Thus, the dimension scales polynomially with $1/\varepsilon$ as $\mathscr{O}((1/\varepsilon)^{D})$. 
Accordingly, the number of qubits required to encode these spaces scales as
$\mathscr{O}\left(
\log D
+
\log(\log(1/\varepsilon))
\right)$ for LHAM, $\mathscr{O}\left(
\log D\log (1/\varepsilon)
\right)$ for Carleman linearization, $\mathscr{O}\left(
D\log (1/\varepsilon)
\right)$ for Koopman-von Neumann approach.
Therefore, the Hilbert space and the qubit complexity in LHAM is exponentially smaller than the ones in Carleman linearization and the Koopman-von Neumann approach. The number of required qubits for LHAM is asymptotically much smaller than those in the other two approaches.
Furthermore, if the solution is approximated with a reduced number of bases as in spectral methods, the complexity of LHAM can be further reduced with a smaller $D$.

The remainder of the paper is organized as follows. In Section \ref{sec:literature}, the relevant work of simulating nonlinear dynamics on quantum computers is reviewed. The details of the proposed LHAM are provided in Section \ref{sec:method} with its error bound theorems. In Section \ref{sec:demonstration}, LHAM and its error scaling are demonstrated with two nonlinear PDEs, including the Burgers' equation and magnetohydrodynamics equations. The paper is concluded in Section \ref{sec:conclusion}. The proofs of theorems related to the LHAM accuracy and the detailed derivations of main results are provided in Appendices \ref{sec:Duhamel}-\ref{sec:MHD}.

\section{\label{sec:literature}Existing Quantum Nonlinear ODE/PDE Solvers}

In the early quantum algorithms to solve nonlinear differential equations, multiple copies of state variables are encoded and computed.  Leyton and Osborne \cite{leyton2008quantum} proposed to encode the original system state variables in the quantum state amplitudes and generate multiple copies of the solution state to evaluate the nonlinear terms including the time derivatives. Because these copies are recursively updated during each time step of Euler method, the required number of initial-state preparations grows exponentially with the number of simulated time steps.
Lloyd et al. \cite{lloyd2020quantum} subsequently proposed a mean-field construction using multiple interacting copies of the solution state, for which the required number of copies scales quadratically with the integration time under suitable stability assumptions.

In the recent methods, nonlinear operators are projected into an infinite-dimensional linear functional space before quantum algorithms are applied to solve linearized systems.
For instance, in the Koopman-von Neumann approach, the original nonlinear dynamics is converted to the Liouville equation, which is a linear evolution equation for a wavefunction defined over the $D$-dimensional continuous state space \cite{joseph2020koopman, dodin2021applications, giannakis2022embedding, jin2022quantum, jin2023time}. The associated probability density satisfies the linear Liouville equation and describes the transport of an ensemble of states under the original nonlinear flow. 
For nonlinear PDEs, spatial discretization first converts the PDE into a $D$-dimensional finite nonlinear ODE system.
The wavefunction is subsequently defined over the $D$ continuous variables rather than over the original physical coordinates. Consequently, if each of these $D$ state-space coordinates is discretized using $N$ grid points, the discretization requires $N^D$ degrees of freedom. This exponential dependence of the discretized Hilbert-space dimension constitutes the curse of dimensionality of discretizations.
In Carleman linearization, the nonlinear differential equation is linearized by introducing additional state variables which are monomials of the original state variables and  the dimension of Hilbert space is enlarged \cite{engel2021linear, liu2021efficient, akiba2023carleman, itani2024quantum, wu2025quantum, brustle2025quantum, gonzalez2025quantum}. 
When truncated at the polynomial degree $N$, the resulting linear system has the dimension $
\sum_{n=1}^{N}\binom{D+n-1}{n}
=
\binom{D+N}{N}-1$.
For a fixed $N$ and $D\gg N$, this dimension scales as
$\mathscr O(D^N/N!)$, which grows polynomially with $D$ for a fixed truncation degree $N$. If $D$ and $N$ instead increase at comparable rates, for example $D=N$, the central-binomial asymptotic of the dimension
$\binom{2D}{D}-1
\sim 4^D/\sqrt{\pi D}$ exhibits exponential growth.
Thus, both linearization methods of Koopman-von Neumann and Carleman linearization suffer from the dimensional blow-up because the Hilbert space dimension grows exponentially.

There have also been efforts to solve nonlinear differential equations using variational quantum algorithms or quantum machine learning \cite{lubasch2020variational, kyriienko2021solving, jaksch2023variational,  su2025quantum}. 
In those algorithms, the solution is obtained by searching for ansatz parameters in a variational manifold or parametric subspace.
Although linearization in general is not necessary in those variational approaches, the variational algorithms are heuristic and prone to optimization difficulties such as barren plateaus and lack of overlap between the solution state and the parametric subspace \cite{mcclean2018barren}. 

\subsection{Homotopy based quantum methods}

Recently, homotopy methods were proposed to solve nonlinear differential equations on quantum computers.
Xue et al. \cite{xue2021quantum} developed an algorithm to solve quadratic ODEs using the homotopy perturbation method. 
The algorithm involved two linearization techniques. The primary one converts a nonlinear PDE into a system of recursive deformation equations, whereas the secondary one embeds the recursive equations into a lifted linear ODE. 
The secondary linearization is to solve the entire linear system at once as homotopy methods exhibit exponential time complexity with the homotopy truncation order when quantum simulation is iteratively applied to each deformation equation, because the homotopy-series solution depends on the earlier ones. 
The secondary linearization was achieved by introducing additional solution fields. 
The resulting linear ODEs are solved step-wise in time with the linear combination of unitaries \cite{berry2017quantum}.
However, the solution of the original nonlinear ODE was obtained with a success probability, which requires additional amplitude amplification.
In addition, similar to Carleman linearization, the secondary linearization increases the dimension of the Hilbert space combinatorially.

Xue et al. \cite{xue2025quantum} later proposed a quantum homotopy analysis method where the secondary linearization involves embedding the recursive linear deformation equations into a linear PDE system defined on a larger spatial domain.
The embedding was performed by introducing auxiliary spatial variables.
With the secondary linearization, the computational complexity of the method increases polynomially with respect to the homotopy truncation order.
However, similar to the previous method \cite{xue2021quantum}, this method experiences the issues of the Hilbert space dimension increasing combinatorially and the query complexity depending on the success probability.

Bharadwaj et al. \cite{bharadwaj2025quantum} developed a similar quantum homotopy algorithm where the deformation equations are embedded in a larger linear system.
The resulting linear system was solved with a time marching compact quantum circuit \cite{bharadwaj2025compact}, which is also implemented with the linear combination of unitaries \cite{childs2012hamiltonian}. 
The linearization of the recursive deformation equations is done by defining the solution with products of lower-order homotopy-series solutions.
Similar to the previous two methods, this approach exhibits issues of the Hilbert space dimension increasing combinatorially and the query complexity depending on the success probability.

Liao \cite{liao2026ham} suggested a framework where the recursive linear PDEs from the homotopy analysis method are solved with Schr\"odingerization \cite{jin2023quantum, jin2024quantum}. However, solving each deformation equation separately with quantum computers requires repetitive state preparations, which creates a major complexity bottleneck due to the no-cloning theorem.

\subsection{Lindbladian approaches to nonunitary quantum simulation}

Lindblad dynamics provides a natural framework for simulating nonunitary evolution on quantum computers. 
Markovian open quantum system dynamics generated by Liouvillians can be  simulated with unitary quantum circuits under suitable locality or sparsity assumptions. For instance,
Kliesch et al. \cite{kliesch2011dissipative} formulated a dissipative quantum Church-Turing theorem and approximated the dissipative dynamics by a Trotter decomposition. 
Childs et al. \cite{3179568.3179569} introduced a sparse Stinespring isometry framework in which the quantum channel generated by each local sparse component is implemented through a Stinespring dilation and ancilla trace-out.
Cleve et al. \cite{cleve_et_al:LIPIcs:2017:7477} used the linear combination of unitaries to decompose the Kraus operators of a quantum channel. The complexity scales nearly linear in the evolution time and polylogarithmic in the inverse precision.

Recently, Li et al. \cite{li_et_al:LIPIcs.ICALP.2023.87} developed a higher-order series expansion approach of Lindbladian quantum channel based on Duhamel's principle. They achieved the similar complexity as the previous work of Cleve et al. \cite{cleve_et_al:LIPIcs:2017:7477} that was based on the first-order approximation. Linear combination of unitaries is utilized to approximate the integration.
Chen et al. \cite{chen2023quantum} developed a compressed incoherent Lindbladian simulation algorithm as a subroutine for quantum Gibbs state preparation. Given a block encoding of the Lindblad operators, the Lindbladian channel can be simulated with the scaling comparable to the prior algorithms \cite{li_et_al:LIPIcs.ICALP.2023.87, cleve_et_al:LIPIcs:2017:7477}. 
Ding et al. \cite{ding2024simulating} proposed a Hamiltonian simulation based approach to construct unitary dynamics on an enlarged Hilbert space and obtain the Lindblad evolution by tracing out ancilla qubits. This method avoids the post selection and can be systematically generalized to higher orders.
Pocrnic et al. \cite{pocrnic2025quantum} proposed an iterative qubitization approach to approximate Lindbladian dynamics with interaction Hamiltonian by which the query complexity scales linearly in time and polylogarithmically in inverse precision.
Shang et al. \cite{shang2025fast} adopted Lindblad dynamics to perform quantum phase estimation and thereby reduce the complexity quadratically.

The Lindblad dynamics simulation was recently adopted to solve general linear ODEs with nonunitary dynamics. Shang et al. \cite{shang2025designing} proposed a nondiagonal density matrix encoding approach, in which the solution of a linear ODE is encoded in an off-diagonal block of a density matrix evolving under a Lindbladian. In this construction, the Hamiltonian portion of the Lindbladian represents the anti-Hermitian component of the ODE generator, while the jump operators represent dissipative components.
Hu and Jin \cite{hu2026amplitude} recently demonstrated that the nondiagonal density matrix encoding and linear combination of Hamiltonian simulation can be unified within a common amplitude-phase structure .

\section{\label{sec:method}The Proposed Lindbladian Homotopy Analysis Method}

In the proposed LHAM, the original nonlinear PDE is converted to a recursive sequence of linear PDEs with the homotopy analysis method.
The recursive nonhomogeneous linear PDEs are reformulated as a higher-dimensional lower block triangular system, which is linear homogeneous autonomous, and solved by performing Lindblad dynamics simulation.
The dissipative part of the nonunitary generator is encoded in the Lindblad jump operators. The resulting Lindblad dynamics is implemented through Stinespring dilation by introducing an environment ancilla qubit and evolving the enlarged system under a dilated Hamiltonian.

\subsection{Linearization}

A nonlinear PDE is defined as
\begin{equation}
\label{eq:general_PDE}
    \frac{\partial}{\partial t}\boldsymbol u(\boldsymbol x,t) = \mathcal M\boldsymbol u(\boldsymbol x,t)+ \mathcal N(\boldsymbol u(\boldsymbol x,t)),
\end{equation}
where $\mathcal M$ is a linear operator
and $\mathcal N$ is a nonlinear operator, both of which include spatial differentiations.
The initial condition of Eq. \eqref{eq:general_PDE} is denoted as $\boldsymbol u(\boldsymbol x,0)= \boldsymbol u_0(\boldsymbol x)$.

In the homotopy analysis method, the nonlinear operator is defined as
\begin{equation}
    \mathcal{N}_{H\!A\!M}(\boldsymbol u(\boldsymbol{x},t))
    :=
    \frac{\partial}{\partial t}\boldsymbol u(\boldsymbol{x},t) - \mathcal M \boldsymbol u(\boldsymbol x,t) - \mathcal N(\boldsymbol u(\boldsymbol x,t)),
    \label{eq:nonlinear_operator_HAM_f}
\end{equation}
which is equivalent to Eq. \eqref{eq:general_PDE} when
$\mathcal{N}_{H\!A\!M}(\boldsymbol{u}(\boldsymbol{x},t))=0$.
The linear operator is chosen as
\begin{equation}
    \mathcal{M}_{H\!A\!M} (\boldsymbol u(\boldsymbol{x},t))
    := \frac{\partial}{\partial t}\boldsymbol u(\boldsymbol{x},t) - \mathcal M\,\boldsymbol u(\boldsymbol{x},t).
    \label{eq:linear_operator_HAM_f}
\end{equation}

\subsubsection{Deformation equations}

A homotopy function is defined as $\boldsymbol \Phi(\boldsymbol{x},t;q)$ with the embedding parameter $q\in[0,1]$. The zeroth-order deformation equation is constructed as
\begin{equation}
    (1-q)\mathcal{M}_{H\!A\!M} \!\big[\boldsymbol \Phi(\boldsymbol{x},t;q) - \boldsymbol u^{(0)}(\boldsymbol{x},t)\big]
    =    q\mu\mathcal{N}_{H\!A\!M} \!\big[\boldsymbol \Phi(\boldsymbol{x},t;q)\big],
    \label{eq:zero_order_deformation_f}
\end{equation}
where $\mu$ is the convergence-control parameter, and $\boldsymbol u^{(0)}(\boldsymbol{x},t)$ is an initial guess satisfying $\boldsymbol u^{(0)}(\boldsymbol{x},0) =\boldsymbol u_0(\boldsymbol{x})$.

The homotopy function $\boldsymbol \Phi(\boldsymbol x, t;q)$ is assumed to be an analytic function which can be represented by a convergent power series with respect to $q$ near $q=0$. This power series, which is also known as the homotopy-Maclaurin series, is defined as
\begin{equation}
    \boldsymbol\Phi(\boldsymbol x,t;q)
    = \boldsymbol u^{(0)}(\boldsymbol x, t) + \sum_{m=1}^{\infty} \boldsymbol u^{(m)}(\boldsymbol x, t)\,q^m,
    \label{eq:homotopy_Maclaurin_series_f}
\end{equation}
where 
$\boldsymbol u^{(m)}(\boldsymbol x, t)
:=
(1/m!)\left.(\partial^m \boldsymbol\Phi/\partial q^m)\right|_{q=0}$.
When $\mu$ in Eq. \eqref{eq:zero_order_deformation_f} is properly chosen,  the series in Eq. \eqref{eq:homotopy_Maclaurin_series_f} converges at $q=1$, where the homotopy-series solution is given as $\boldsymbol u(\boldsymbol x, t)= \boldsymbol\Phi(\boldsymbol x, t;1)$.

Since $\mathcal{M}_{H\!A\!M}$ in Eq. \eqref{eq:linear_operator_HAM_f} is independent of $q$, it commutes with $\partial/\partial q$. When Eq. \eqref{eq:zero_order_deformation_f} is differentiated $m$ times with respect to $q$ and evaluated at $q=0$, the $m$th-order deformation equation is constructed as
\begin{equation}
    \mathcal{M}_{H\!A\!M}\!\big[\boldsymbol u^{(m)}(\boldsymbol x, t) - \boldsymbol u^{(m-1)}(\boldsymbol x, t)\big]
    =
    \mu\,\boldsymbol R^{(m-1)}(\boldsymbol x, t)
    \label{eq:kth_deformation_coeff_f}
\end{equation}
for $m \ge 2$, where
\begin{equation}
    \boldsymbol R^{(m-1)}(\boldsymbol x, t)
    :=
    \frac{1}{(m-1)!}
    \left.
    \frac{\partial^{m-1}}{\partial q^{m-1}}
    \mathcal{N}_{H\!A\!M}\!\big[\boldsymbol\Phi(\boldsymbol x, t;q)\big]
    \right|_{q=0}   ,
    \label{eq:R_coeff_def_f}
\end{equation}
and is constructed as
\begin{equation}
    \mathcal{M}_{H\!A\!M}\!\big[\boldsymbol u^{(1)}(\boldsymbol x, t)\big]
    =
    \mu\,\boldsymbol R^{(0)}(\boldsymbol x, t)
    \label{eq:1th_deformation_coeff_f}
\end{equation}
for $m=1$.
Eqs. \eqref{eq:kth_deformation_coeff_f} and \eqref{eq:1th_deformation_coeff_f}  are combined into a system of nonhomogeneous linear PDEs for the unknown $\boldsymbol u^{(m)}(\boldsymbol x,t)$.
For higher-order deformation equations where $m \ge 1$, the initial conditions are set to $\boldsymbol u^{(m)}(\boldsymbol{x},0) = \boldsymbol 0$ so that $\sum_{m\ge 0} \boldsymbol u^{(m)}(\boldsymbol{x},0)=\boldsymbol u_0(\boldsymbol{x})$. 
That is,
\begin{equation}
\label{eq:init_conditions_f}
    \boldsymbol u^{(m)}(\boldsymbol{x},0)=
    \begin{cases}
    \boldsymbol u_0(\boldsymbol{x}) & \text{for } m=0 \\
    \boldsymbol 0 & \text{for } m\ge1 .
    \end{cases}
\end{equation}

By substituting Eq. \eqref{eq:linear_operator_HAM_f} into Eqs. \eqref{eq:kth_deformation_coeff_f} and \eqref{eq:1th_deformation_coeff_f}, the $m$th-order deformation equation is rewritten as
\begin{equation}
    \frac{\partial}{\partial t}\boldsymbol u^{(m)}(\boldsymbol x,t)= \mathcal M\,\boldsymbol u^{(m)}(\boldsymbol x,t)+\boldsymbol f^{(m-1)}(\boldsymbol x,t),
    \label{eq:forced_linear_cm_f}
\end{equation}
with the initial linear homogeneous differential equation $
    \partial\boldsymbol {u}^{(0)}(\boldsymbol x,t)/\partial t
    = \mathcal M\,\boldsymbol u^{(0)}(\boldsymbol x,t)$.
The forcing terms are
$\boldsymbol f^{(0)}(\boldsymbol x,t) := \mu\,\boldsymbol R^{(0)}(\boldsymbol x,t)$ for $m=1$, and  $
    \boldsymbol f^{(m-1)}(\boldsymbol x,t) :=  \partial\boldsymbol u^{(m-1)}/\partial t - \mathcal M\,\boldsymbol u^{(m-1)} + \mu\,\boldsymbol R^{(m-1)}$
for $m\ge2$. 
The $(m-1)$th-order homotopy residual $\boldsymbol{R}^{(m-1)}$ is obtained by substituting Eq. \eqref{eq:nonlinear_operator_HAM_f} into Eq. \eqref{eq:R_coeff_def_f} such that $
    \boldsymbol R^{(m-1)}(\boldsymbol x,t)
    =
    \partial\boldsymbol u^{(m-1)}(\boldsymbol x,t)/\partial t - \mathcal M\,\boldsymbol u^{(m-1)}(\boldsymbol x,t)
    - \boldsymbol S^{(m-1)}(\boldsymbol x,t)$, 
where the nonlinear source term is obtained as $
    \boldsymbol S^{(m-1)}(\boldsymbol x,t)
    :=
    (1/(m-1)!)\left.(\partial^{m-1}\mathcal N(\boldsymbol\Phi(\boldsymbol x,t;q))/\partial q^{m-1})
    \right|_{q=0}$.
As a result, $\boldsymbol f^{(m-1)}$ in the $m$th-order deformation equation can be calculated from the solutions $\{\boldsymbol u^{(0)},\ldots,\boldsymbol u^{(m-1)}\}$ of the deformation equations of previous orders based on $\boldsymbol S^{(m-1)}(\boldsymbol x,t)$ and $\boldsymbol R^{(m-1)}(\boldsymbol x,t)$.
In the traditional homotopy analysis method, the solutions of the deformation equations are obtained recursively from the lower order to higher order.
In the proposed LHAM, all deformation equations are solved simultaneously instead, which will be described in Sec. \ref{sec:lifting}.

The analytic solution of Eq. \eqref{eq:forced_linear_cm_f} is 
\begin{equation}
    \boldsymbol u^{(m)}(\boldsymbol x,t)
    =
    \int_0^t e^{\mathcal M(t-\tau)}\,\boldsymbol f^{(m-1)}(\boldsymbol x,\tau)\,d\tau.
    \label{eq:burgers_duhamel_f}
\end{equation}
The derivation of Eq. (\ref{eq:burgers_duhamel_f}) is based on the the Duhamel's principle, given in Appendix \ref{sec:Duhamel}.

\subsubsection{Lifted linearization of the recursive nonhomogeneous linear PDEs}\label{sec:lifting}

To obtain a closed linear homogeneous autonomous system, we propose to represent each forcing term as
\begin{equation}
    \boldsymbol f^{(m-1)}(\boldsymbol x,t) = \mathcal V^{(m)} \boldsymbol z^{(m)}(\boldsymbol x,t),
    \label{eq:single_channel_force_f}
\end{equation}
where $\mathcal V^{(m)}$ is a feed-forward coupling operator to map between $\boldsymbol f^{(m-1)}$ and a new auxiliary vector $\boldsymbol z^{(m)}$.
The auxiliary vector $\boldsymbol z^{(m)}(\boldsymbol x,t)=e^{-\mathbf \Lambda^{(m)} t}\boldsymbol z^{(m)}(\boldsymbol x,0)$ is introduced such that its dynamics is described by 
\begin{equation}
\label{eq:zm_evolution}
    \frac{\partial}{\partial t} \boldsymbol z^{(m)}(\boldsymbol x,t)=-\mathbf \Lambda^{(m)} \boldsymbol z^{(m)}(\boldsymbol x,t),
\end{equation}
where $\mathbf \Lambda^{(m)}$ is an auxiliary evolution operator that generates the decaying magnitudes of $\boldsymbol{z}^{(m)}(\boldsymbol{x},t)$.

The computation is based on the lifted state
\begin{equation}
    \boldsymbol y(t)
    =
    \begin{bmatrix}
        \boldsymbol u^{(0)}(\boldsymbol x,t)\\
        \boldsymbol z^{(1)}(\boldsymbol x,t)\\
        \boldsymbol u^{(1)}(\boldsymbol x,t)\\
        \boldsymbol z^{(2)}(\boldsymbol x,t)\\
        \vdots\\
        \boldsymbol z^{(\tilde{m})}(\boldsymbol x,t)\\
        \boldsymbol u^{(\tilde{m})}(\boldsymbol x,t)
    \end{bmatrix},
    \label{eq:mu_def_f}
\end{equation}
where $\tilde{m}$ is the homotopy truncation order.
With the definition in Eq. \eqref{eq:mu_def_f}, instead of the traditional recursive process to solve deformation equations in increasing orders, all deformation equations are combined into one linear homogeneous autonomous system
\begin{equation}
    \frac{d}{dt}{\boldsymbol y}(t)= \mathcal A\,\boldsymbol y(t),
    \label{eq:augmented_linear_system_f}
\end{equation}
where
\begin{equation}
    \mathcal A=
    \begin{bmatrix}
        \mathcal M & 0 & 0 & 0 & \cdots & 0 & 0\\
        0 & -\mathbf \Lambda^{(1)} & 0 & 0 & \cdots & 0 & 0\\
        0 & \mathcal V^{(1)} & \mathcal M & 0 & \cdots & 0 & 0\\
        0 & 0 & 0 & -\mathbf \Lambda^{(2)} & \cdots & 0 & 0\\
        0 & 0 & 0 & \mathcal V^{(2)} & \cdots & 0 & 0\\
        \vdots & \vdots & \vdots & \vdots & \ddots & -\mathbf \Lambda^{(\tilde{m})} & 0\\
        0 & 0 & 0 & 0 & \cdots & \mathcal V^{(\tilde{m})} & \mathcal M
    \end{bmatrix}
    \label{eq:A_block_form_f}
\end{equation}
is a lower block triangular operator.
The diagonal blocks $\mathcal M$'s evolve state variable $\boldsymbol u^{(m)}(\boldsymbol x, t)$ for each $m$ under the same linear operator.    
The behavior of the linear homogeneous autonomous system in Eq. \eqref{eq:augmented_linear_system_f} is a generic nonunitary dynamics described as
\begin{equation}
\label{eq:solution_f}
    \boldsymbol y(t)=e^{\mathcal At}\boldsymbol y(0).
\end{equation}
Eq. \eqref{eq:augmented_linear_system_f} can be solved on a quantum computer by discretizing and encoding $\boldsymbol y(t)$ in Eq. \eqref{eq:mu_def_f}. Alternatively, spectral methods can be used to solve Eq. \eqref{eq:augmented_linear_system_f} by encoding the coefficients of bases to further reduce the number of qubits.

\subsubsection{Closure of the auxiliary variables under resonant forcing}

If $-\sigma_k$ is an eigenvalue of $\mathcal M$, then the time evolution in the corresponding eigenspace of $\mathcal M$ is represented as $e^{-\sigma_k t}$. 
The time evolution exhibits resonant forcing when the decay rate $-\sigma_k$ is equivalent to $\lambda$, which is an eigenvalue of the auxiliary evolution operator $\boldsymbol{\Lambda}^{(m)}$. 
In this case, a particular solution is generated by a resonant forcing basis function $
    \boldsymbol f_{\lambda,p}(\boldsymbol x,t)
    =
    \boldsymbol b_{\lambda,p}(\boldsymbol x)
    t^p e^{-\lambda t}$,
where $\boldsymbol b_{\lambda,p}$ is time-independent coefficient function and $p$ is a nonnegative integer. The particular solution contains a term proportional to $t^{p+1}e^{-\lambda t}$, where the order of auxiliary basis function is increased by one. Therefore, the time derivative of every auxiliary variable can be represented as a linear combination of variables within the same finite-dimensional auxiliary system.
By linearity of the deformation equation, the response to each resonant forcing basis function can be determined separately and then added together. From Eq.~\eqref{eq:burgers_duhamel_f}, the response corresponding to the initial condition in Eq. \eqref{eq:init_conditions_f} is
\begin{equation}
    \boldsymbol u_{\lambda,p}(\boldsymbol x,t)
    =
    \int_0^t
    e^{\mathcal M(t-\tau)}
    \boldsymbol b_{\lambda,p}(\boldsymbol x)
    \tau^p e^{-\lambda\tau}\,d\tau.
    \label{eq:operator_forcing_response}
\end{equation}

The response in Eq. \eqref{eq:operator_forcing_response} can be projected onto the eigenspace of $\mathcal M$ associated with the eigenvalue $-\sigma_k$ with a projection operator $\mathcal P_k$ such that $
    \mathcal P_k e^{t\mathcal M}
    =
    e^{-\sigma_k t}\mathcal P_k$.
Applying $\mathcal P_k$ to
Eq.~\eqref{eq:operator_forcing_response} yields
\begin{equation}
    \mathcal P_k\boldsymbol u_{\lambda,p}(\boldsymbol x,t)
    =
    e^{-\sigma_k t}
    \mathcal P_k\boldsymbol b_{\lambda,p}(\boldsymbol x)
    \int_0^t
    \tau^p e^{(\sigma_k-\lambda)\tau}\,d\tau.
    \label{eq:projected_forcing_response}
\end{equation}
If $\lambda=\sigma_k$, then Eq. \eqref{eq:projected_forcing_response} simplifies to
\begin{equation}
    \mathcal P_k\boldsymbol u_{\lambda,p}(\boldsymbol x,t)
    =
    \frac{t^{p+1}}{p+1}e^{-\lambda t}
    \mathcal P_k\boldsymbol b_{\lambda,p}(\boldsymbol x).
    \label{eq:operator_resonant_response}
\end{equation}
In contrast to the nonresonant case (i.e., $p=0$ and $\lambda\neq\sigma_k$), where 
\begin{equation}
\label{eq:operator_nonresonant_response}
    \mathcal P_k\boldsymbol u_{\lambda,0}(\boldsymbol x,t)
    =
    \frac{e^{-\lambda t}-e^{-\sigma_k t}}
    {\sigma_k-\lambda}
    \mathcal P_k\boldsymbol b_{\lambda,0}(\boldsymbol x)
\end{equation}
consists of exponential terms with respect to time,
the resonant case involves responses which are both exponentially and monomially dependent on time. 
It is also observed in Eq.~\eqref{eq:operator_resonant_response} that the mononomial degree of the resonant forcing response increases by one. 
The forcing response terms at higher deformation orders are formulated similarly to Eq. \eqref{eq:operator_resonant_response}. The auxiliary basis is expanded to include functions which are exponentailly and monomially dependent on time.

It is assumed that, through the retained homotopy order, the forcing admits a finite representation of the form
$
    \boldsymbol f^{(m-1)}(\boldsymbol x,t)
    =
    \sum_{\lambda\in\mathbb A}
    \sum_{p=0}^{r_\lambda}
    \boldsymbol v_{\lambda,p}^{(m)}(\boldsymbol x)
    z_{\lambda,p}(t)$,
where $z_{\lambda,p}(t)=t^p e^{-\lambda t}$ is the temporal basis functions as scalar auxiliary variables, whose autonomous evolution generates the time dependence of the forcing. $\mathbb A$ is the set of all distinct exponential rates required by the retained deformation equations, including both resonant and nonresonant rates.
For each $\lambda$, $r_\lambda$ is the largest monomial degree required in this forcing representation.
The coefficient functions $\boldsymbol v_{\lambda,p}^{(m)}(\boldsymbol{x})$ are time-independent and are zero when the corresponding temporal basis function does not contribute to the forcing at a particular deformation order $m$.

The auxiliary variables satisfy
\begin{align}
\label{eq:auxiliary_polynomial_closure}
\begin{split}
    \frac{d}{dt}z_{\lambda,0}(t)
    &=
    -\lambda z_{\lambda,0}(t),
    \\
    \frac{d}{dt}z_{\lambda,p}(t)
    &=
    -\lambda z_{\lambda,p}(t)
    +p z_{\lambda,p-1}(t),
\end{split}
\end{align}
where $p=1,\ldots,r_\lambda$.
Thus, differentiation maps each auxiliary variable into the span of variables with the same exponential rate.
Including all degrees from $0$ to $r_\lambda$ therefore makes this span invariant under differentiation.

For each $\lambda\in\mathbb A$, the column vector is defined as
\begin{equation}
    \boldsymbol z_\lambda(t)
    =
    \begin{pmatrix}
        z_{\lambda,0}(t),
        z_{\lambda,1}(t),
        \dots,
        z_{\lambda,r_\lambda}(t)
    \end{pmatrix}.
\end{equation}
Eq.~\eqref{eq:auxiliary_polynomial_closure} becomes
\begin{equation}
    \frac{d}{dt}\boldsymbol z_\lambda(t)
    =
    -\Lambda_\lambda\boldsymbol z_\lambda(t),
    \qquad
    \boldsymbol z_\lambda(0)
    =
    \begin{pmatrix}
        1, 0, \dots, 0
    \end{pmatrix},
\end{equation}
where
\begin{equation}
    \Lambda_\lambda
    =
    \begin{bmatrix}
        \lambda & 0 & 0 & \cdots & 0\\
        -1 & \lambda & 0 & \cdots & 0\\
        0 & -2 & \lambda & \cdots & 0\\
        \vdots & \vdots & \vdots & \ddots & \vdots\\
        0 & 0 & \cdots & -r_\lambda & \lambda
    \end{bmatrix}
    \label{eq:auxiliary_rate_block}
\end{equation}
is a monomial-exponential auxiliary variable generator.
When $r_\lambda=0$, this block reduces to the scalar
$\lambda$. When $r_\lambda>0$, its off-diagonal entries generate the monomial factors multiplying
$e^{-\lambda t}$.

By collecting the auxiliary variables, the auxiliary vector $
\boldsymbol z(t)
=
\bigoplus_{\lambda\in\mathbb A}\boldsymbol z_\lambda(t)$ is obtained. The auxiliary evolution operator $
\mathbf\Lambda
=
\bigoplus_{\lambda\in\mathbb A}\Lambda_\lambda$ is block diagonal because differentiation preserves each exponential rate.
The enlarged auxiliary vector satisfies $   d\boldsymbol z(t)/dt
    =
    -\mathbf\Lambda\boldsymbol z(t)$ and $\boldsymbol z(t)
    =
    e^{-\mathbf\Lambda t}\boldsymbol z(0)$.    
The same auxiliary variables can be shared across
deformation orders, with their contributions specified
by order-dependent coupling operators.

$\mathcal V^{(m)}$ is denoted as the combined coupling
operator defined by
\begin{equation}
    \bigl(\mathcal V^{(m)}\boldsymbol z(t)\bigr)(\boldsymbol x)
    =
    \sum_{\lambda\in\mathbb A}
    \sum_{p=0}^{r_\lambda}
    \boldsymbol v_{\lambda,p}^{(m)}(\boldsymbol x)
    z_{\lambda,p}(t).
\end{equation}
The deformation equation retains the form of Eq.~\eqref{eq:forced_linear_cm_f}, with the forcing term $\boldsymbol{f}^{(m-1)}$ represented by $\mathcal V^{(m)}\boldsymbol z(t)$.
To avoid storing repeated copies of the same auxiliary
variables, the order-specific auxiliary vectors are
consolidated into the shared vector $\boldsymbol z(t)$.
The augmented evolution remains a linear homogeneous
autonomous system by replacing 
Eq. \eqref{eq:A_block_form_f} with
\begin{equation}
    \mathcal A
    =
    \begin{bmatrix}
        -\mathbf\Lambda & 0 & 0 & \cdots & 0\\
        0 & \mathcal M & 0 & \cdots & 0\\
        \mathcal V^{(1)} & 0 & \mathcal M & \cdots & 0\\
        \vdots & \vdots & \vdots & \ddots & \vdots\\
        \mathcal V^{(\tilde m)} & 0 & 0 & \cdots & \mathcal M\\
    \end{bmatrix},
    \label{eq:A_block_form_f_res}
\end{equation}
where the shared auxiliary vector is coupled to each deformation order.
Eq. \eqref{eq:mu_def_f} is replaced by 
\begin{equation}
    \boldsymbol y(t)
    =
    \begin{bmatrix}
        \boldsymbol z(t) \\
        \boldsymbol u^{(0)}(\cdot,t)\\
        \boldsymbol u^{(1)}(\cdot,t)\\
        \vdots\\
        \boldsymbol u^{(\tilde m)}(\cdot,t)
    \end{bmatrix},
    \label{eq:mu_def_f_res}
\end{equation}
where $\boldsymbol u^{(m)}(\cdot,t)$ denotes the spatial
function at time $t$.

With a shared temporal basis, the number of scalar auxiliary variables is $
    d_{\mathrm{A}}
    =
    \sum_{\lambda\in\mathbb A}(r_\lambda+1)$.
For the subset $\mathbb S_{\mathrm{R}}\subseteq\mathbb A$ of rates participating in resonance, the corresponding auxiliary variables has dimension
$
    d_{R}
    =
    \sum_{\lambda\in\mathbb S_{R}}
    (r_\lambda+1)$.
If the pure-exponential variables $z_{\lambda,0}$ are already included, augmenting these blocks requires
$
    \hat d_{\mathrm R}
    =
    \sum_{\lambda\in\mathbb S_{R}}r_\lambda$ additional variables.
These counts depend on the distinct exponential rates and the largest monomial degree required for each rate.
The same temporal auxiliary variable can be shared across deformation orders through the spatial coupling functions $\boldsymbol v_{\lambda,p}^{(m)}(\boldsymbol x)$ defining the order-dependent coupling operators $\mathcal V^{(m)}$.
Auxiliary variables associated with monomially weighted nonresonant rates, when required, are also included in $d_{\mathrm A}$.

After spatial discretization, $D$ is denoted as the dimension of each physical deformation vector.
Because the shared auxiliary vector is stored only once, the total lifted-state dimension is
$ 
    d_{T}
    =
    D(\tilde m+1)
    +
    \sum_{\lambda\in\mathbb A}(r_\lambda+1)$.
If the number of distinct required rates and their maximum monomial degree satisfy
$
    |\mathbb A|=\mathscr O(D)$ and
    $
    \max_{\lambda\in\mathbb A}r_\lambda
    =
    \mathscr O(\tilde m)$,
then
$
    d_T
    =
    \mathscr O\!\left(D(\tilde m+1)\right)
    =
    \mathscr O(D\tilde m)$, where
    $\tilde m\geq 1$.
Under these assumptions, incorporating resonance preserves the lifted-state dimension scaling.

\subsection{Lindblad dynamics}

To simulate the nonunitary evolution in Eq. \eqref{eq:solution_f}, the linear dissipative operator is encoded in the jump operator in the Lindblad master equation \cite{shang2025designing}. 
The Lindblad master equation can be solved based on the Stinespring dilation. The evolved state in Eq. \eqref{eq:solution_f} is also embedded into the upper-right block of the updated density matrix for each time step. 

The generator $\mathcal A$ is decomposed into its Hermitian and anti-Hermitian components, as $
    \mathcal A_1= -(\mathcal A+\mathcal A^\dagger)/2$ and
    $\mathcal A_2= -(\mathcal A-\mathcal A^\dagger)/(2i)$,
respectively, where $\mathcal A_1$ generates the dissipative component of the dynamics with a norm contraction, whereas $\mathcal A_2$ generates the unitary rotations.

For the lifted operator $\mathcal A$, the Hermitian component $\mathcal A_1$ is not necessarily positive semidefinite. To restore the semi-dissipative structure, we introduce a scalar shift of the generator so that $\hat{\mathcal{A}} = \mathcal A - \gamma \mathcal I$, with $\gamma \ge - \mathrm{inf}\,\mathrm{spec}(\mathcal A_1)$,
where $\mathrm{spec}(\mathcal A_1)$ is the spectrum of $\mathcal A_1$. As a result, we obtain $\hat{\mathcal{A}}_1 = \mathcal A_1 + \gamma \mathcal I \succeq 0$, and denote $\hat{\mathcal{A}}_2 = \mathcal A_2$.
The solution of the shifted ODE is related to the original solution by $\hat{\boldsymbol y}(t) = e^{-\gamma t}\boldsymbol y(t) $.

We embed the problem in an enlarged system space by introducing an ancilla basis $\{|0\rangle,|1\rangle\}$. The system Hamiltonian and jump operator are defined as \cite{shang2025designing, hu2026amplitude} 
\begin{equation}
\label{eq:hamiltonian_jumpop}
    \mathcal H=
    \begin{bmatrix}
        \hat{\mathcal{A}}_2 & 0\\
        0 & 0
    \end{bmatrix}
    \; \text{and} \;
    \mathcal F=
    \begin{bmatrix}
        \sqrt{2\hat{\mathcal{A}}_1} & 0\\
        0 & 0
    \end{bmatrix},
\end{equation}
respectively. The density operator $\rho$ evolves according to the Lindblad master equation 
\begin{equation}
\label{eq:lindbladian}
    \frac{d\rho}{dt}
    =
    \mathcal L(\rho)
    =
    -i[\mathcal H,\rho]
    +
    \sum_{k=1}^{K}
    \left(
        \mathcal F_k\rho \mathcal F_k^\dagger
        -
        \frac12\{\mathcal F_k^\dagger \mathcal F_k,\rho\}
    \right),
\end{equation}
where $\mathcal{L}$ is the Lindbladian superoperator. In this paper, the jump operator $\mathcal F$ in Eq. \eqref{eq:hamiltonian_jumpop} is defined for $K=1$.
The upper-right block is 
\begin{equation}
\label{eq:rho_01}
\rho_{01}(t) =
(\bra{0}\otimes \mathcal I)\rho(t)(\ket{1}\otimes \mathcal I).
\end{equation}

For the initial state
$\rho(0)=\ket{+}\bra{+}\otimes \ket{\boldsymbol y(0)}\bra{\boldsymbol y(0)}$,
where $|\boldsymbol y(0)\rangle$ is normalized,
the upper-right block is
$\rho_{01}(0)
=
(1/2)
\ket{\boldsymbol y(0)}\bra{\boldsymbol y(0)}$.
In the Lindblad dynamics, the upper-right block of $\rho$ evolves as
\begin{equation}
\label{eq:rho01T_f}
    \rho_{01}(t)
    =
    \frac{1}{2}
    |\hat{\boldsymbol y}(t)\rangle
    \langle \boldsymbol y(0) |,
\end{equation}
It is seen in Eq. \eqref{eq:rho01T_f} that the solution of the original system is encoded in the upper-right block of the density operator.

The evolution of the system is equivalent to $d\rho/dt = \mathcal L(\rho)$.
The one-step quantum channel is completely positive trace preserving map, which can be written with Kraus operators as
\begin{equation}
\label{eq:cptp}
    \mathcal E_{\Delta t}(\rho(t)) := e^{\mathcal L \Delta t}(\rho(t)) = \rho(t+\Delta t) = \sum_{s=0}^S \mathcal K_s \rho(t) \mathcal K_s^\dag,
\end{equation}
and applied repeatedly for $\tilde n$ steps. 
The Kraus operators are defined as $\mathcal K_0 = \mathcal I + \mathcal F_0\Delta t$, $\mathcal K_1 = \mathcal F\sqrt{\Delta t}$ for $S=1$, and  $\mathcal K_0 = \mathcal I + \mathcal F_0\Delta t + (1/2)\mathcal F_0^2 \Delta t^2$, $\mathcal K_1 = -i\sqrt{\Delta t}(\mathcal F + (\Delta t/2)(\mathcal F\mathcal F_0 + \mathcal F_0\mathcal F))$, $\mathcal K_2 = -i(\sqrt{2}\Delta t/2)\mathcal F^2$
for $S=2$, where $\mathcal F_0 = -i\mathcal H -(1/2)\mathcal F^\dag \mathcal F$ \cite{ding2024simulating}.

The total evolution time is $T = \tilde n\Delta t$.
Every completely positive trace preserving map in Eq. (\ref{eq:cptp}) can be implemented by Stinespring dilation where environment ancilla qubits $\ket{0}_E$ are introduced, as
\begin{equation}
\label{eq:stinespring_dilation}
     \mathcal E_{\Delta t}(\rho_n)
     = \rho_{n+1}
     =\mathrm{Tr}_E \bigl(\mathcal U(\ket{0}\bra{0}_E \otimes \rho_n)\mathcal U^\dagger \bigr),
\end{equation}
where $\mathcal U$ is unitary operator evolving the system-environment state and $\rho$ is the evolved system density operator. 
The unitary operator $\mathcal U = e^{-i\sqrt{\Delta t}\widetilde{\mathcal H}}$ is constructed with the dilated Hamiltonian \cite{ding2024simulating}
\begin{equation}
\label{eq:ding_dilated_hamiltonian}
    \widetilde{\mathcal H}
    =
    |0\rangle\langle 0|\otimes \mathcal H_0
    +
    \sum_{s=1}^{S}
    \left(
        |s\rangle\langle 0|\otimes \mathcal H_s
        +
        |0\rangle\langle s|\otimes \mathcal H_s^\dagger
    \right),
\end{equation}
where the $\mathcal H_s$ blocks are determined by expanding $\mathcal U$ in powers of $\sqrt{\Delta t}$ and matching, order by order, the
first-column blocks $
    {}_E\langle s|\mathcal U|0\rangle_E$
with the corresponding Kraus operators $\mathcal K_s$ of the order-$S$ approximation to the quantum channel in Eq.~\eqref{eq:cptp}.
For $S=1$, matching the Taylor expansion through the terms required for a local error of order $\mathscr O(\Delta t^2)$ gives 
$\mathcal H_0 = \Delta t^{1/2} \mathcal H$, $\mathcal H_1 = \mathcal F$. For $S=2$, matching the Taylor expansion through the terms required for a local error of order $\mathscr O(\Delta t^3)$ gives $\mathcal H_0 = \Delta t^{1/2}\mathcal H -(\Delta t^{3/2}/12)\{\mathcal H,\mathcal F^\dag \mathcal F\}$, $\mathcal H_1 = \mathcal F + \Delta t\big((1/2)\mathcal F_0\mathcal F + (1/2)\mathcal F\mathcal F_0 + (1/6)\mathcal F\mathcal F^\dag \mathcal F + (i/2)\mathcal F\mathcal H \big)$, $\mathcal H_2 = (\Delta t^{1/2}/\sqrt{2})\mathcal F^2 + \Delta t^{3/2} \big((1/6\sqrt{2})\mathcal F^2\mathcal F^\dag \mathcal F + (i/2\sqrt{2})\mathcal F^2\mathcal H \big)$.

The partial trace over the environment can be implemented efficiently by the mid-circuit-measure-and-reset strategy.
After the final step, 
the upper-right block $\rho_{01}(T)$ is extracted by measuring $\rho$ on the ancilla quit basis as in Eq. \eqref{eq:rho_01} and each component of the solution $\boldsymbol y(T)$ is recovered by measurements of Hermitian observables.

The Lindblad dynamics simulation results in $\rho(T) \in \mathbb{H}_a \otimes \mathbb{H}_s$, where $\mathbb{H}_a$ and $\mathbb{H}_s$ are the Hilbert spaces corresponding to the ancilla and the original system, respectively. The upper-right block operator $\rho_{01}$ in Eq. \eqref{eq:rho_01} is equivalent to $
    \rho_{01}=(1/2)\bigl[\mathrm{Tr}_a((X\otimes \mathcal I)\,\rho) - i\,\mathrm{Tr}_a((Y\otimes \mathcal I)\,\rho) \bigr]$,
where $X=\ket{0}\bra{1}+\ket{1}\bra{0}$ and $Y= -i\ket{0}\bra{1}+i\ket{1}\bra{0}$. 
Each component of $\ket{\hat{\boldsymbol y}(t)}=\sum_j \hat{y}_j(t)\ket{j}$ is extracted as $\hat{y}_j(t) = 2\bra{j}\rho_{01}\ket{\boldsymbol{y}(0)} =2\,\mathrm{Tr}(\ket{\boldsymbol{y}(0)}\bra{j}\rho_{01})$.
Equivalently, each component is measured as
\begin{align}\label{eq:y_tilde_Tr}
\begin{split}
    \hat{y}_j(t) 
    &= \frac{1}{2}\bigl[\mathrm{Tr}\bigl((X\otimes \mathcal O_j^{(R)})\,\rho\bigr) + i\,\mathrm{Tr}\bigl((X\otimes \mathcal O_j^{(I)})\,\rho\bigr) \\
    &- i\,\mathrm{Tr}\bigl((Y\otimes \mathcal O_j^{(R)})\,\rho\bigr)
    + \mathrm{Tr}\bigl((Y\otimes \mathcal O_j^{(I)})\,\rho\bigr)\bigr],
\end{split}
\end{align}
where the measurable Hermitian observables are defined as
\begin{equation}
\label{eq:observable}
\begin{split}
    \mathcal O_j^{(R)} &:= \ket{\boldsymbol{y}(0)}\bra{j} + \ket{j}\bra{\boldsymbol{y}(0)}, \\
    \mathcal O_j^{(I)} &:= -i\,\ket{\boldsymbol{y}(0)}\bra{j} + i\,\ket{j}\bra{\boldsymbol{y}(0)}.
\end{split}
\end{equation}
Then, $\mathrm{Re}(\bra{j}\rho_{01}\ket{\boldsymbol{y}(0)})$ and $\mathrm{Im}(\bra{j}\rho_{01}\ket{\boldsymbol{y}(0)})$ are obtained separately from the expectation values in Eq. \eqref{eq:y_tilde_Tr} involving the Hermitian operators in Eq. \eqref{eq:observable}.
Finally, each component $y_j(t)$ of $\boldsymbol{y}(t)$ is recovered as $y_j(t) = e^{\gamma t}\hat{y}_j(t)$.
Therefore, the solution vector is recovered by measuring Hermitian observables, which does not involve post-selection success probability.

\subsection{Functional expansion}
To reduce the number of qubits, the functional expansion  \cite{sul2024quantum, sul2025generic} can be used.
A scalar field solution $ u(\boldsymbol x,t)$ is approximated as
\begin{equation}
\label{eq:Galerkin}
    u(\boldsymbol x,t) \approx \tilde{u}(\boldsymbol x,t) := \sum_{j=0}^{J-1} c_{j}(t)\varphi_{j}(\boldsymbol x),
\end{equation}
where $J$ is the number of basis functions for the truncated $\tilde{u}$, $\{\varphi_j\}_{j\ge0}$ is a set of orthonormal basis functions defined on $\mathbb H_J := \mathrm{span}\{\varphi_0,\ldots,\varphi_{J-1}\}$, and the coefficients are calculated based on the $L^2$ inner product $
    c_j(t) = \langle \varphi_j, u(\cdot,t)\rangle
    := \int_{-\infty}^{\infty} \varphi_j^*(\boldsymbol x)u(\boldsymbol x,t)\,d\boldsymbol{x}$.
Since $\varphi_j(\boldsymbol x)$ is time-independent, the time derivative commutes with the
inner product such as
$dc_j(t)/dt 
    = d\langle\varphi_j,u(\cdot,t)\rangle /dt
    = \langle\varphi_j,\partial u(\cdot,t)/\partial t \rangle$.
The linear term in Eq. \eqref{eq:general_PDE} is projected onto the basis such that
    $\langle \varphi_k, \mathcal M \tilde{u}\rangle
    = \sum_{j=0}^{J-1} M_{kj}\,c_j(t)$
    for $0\le k\le J-1$
where 
    $M_{kj} = \langle \varphi_k, \mathcal M\varphi_j\rangle$
is the element of matrix $M =[M_{kj}]\in\mathbb{C}^{J \times J} $ .
The $j$th component in the expansion of the nonlinear term in Eq. \eqref{eq:general_PDE} is $
    \langle \varphi_j, \mathcal N(\tilde{u})\rangle
    =
    N_j(\boldsymbol{c}(t))$,
which is a nonlinear function of $\boldsymbol{c}(t)=(c_{0}(t),\ldots,c_{J-1}(t)) \in\mathbb C^J$.

The coefficients after the functional expansion follow ODEs
\begin{equation}
    \frac{d}{dt}{\boldsymbol c}(t)
    = M\,\boldsymbol c(t) + \,\boldsymbol{N}(\boldsymbol c(t)),
    \label{eq:ode_coeff}
\end{equation}
where $\boldsymbol{N}(\boldsymbol{c}(t)) = (N_0(\boldsymbol{c}(t)), \ldots, N_{J-1}(\boldsymbol{c}(t)))$.
The initial condition is $\boldsymbol c(0)=\boldsymbol c_{0}$ with
$\boldsymbol c_0 = [c_{0,j}] \in \mathbb{C}^J$ where $c_{0,j} = \langle \varphi_j,u_0\rangle$.
The homotopy-Maclaurin series is expanded as Eq. \eqref{eq:homotopy_Maclaurin_series_f} by replacing the variables $\boldsymbol{u}(\boldsymbol{x},t)$ by $\boldsymbol{c}(t)$. 
The $m$th-oder deformation equations are rewritten as
\begin{equation}
    \frac{d}{dt}\boldsymbol c^{(m)}(t)= M\,\boldsymbol c^{(m)}(t)+\boldsymbol f^{(m-1)}(t)
    \label{eq:forced_linear_cm}
\end{equation}
with the initial condition $\boldsymbol c^{(m)}(0)=\boldsymbol 0$. The initial linear homogeneous differential equation is
\begin{equation}
    \frac{d}{dt}\boldsymbol {c}^{(0)}(t)
    =M\,\boldsymbol c^{(0)}(t)
    \label{eq:c0_hom}
\end{equation}
with $\boldsymbol{c}^{(0)}(0)=\boldsymbol{c}_0$.

To obtain a closed linear system, each forcing term is represented in the form of Eq. \eqref{eq:single_channel_force_f}.
The nonlinear forcing term can be represented with finite auxiliary variables which evolve as in Eq. \eqref{eq:zm_evolution}.
The auxiliary variables are grouped into a vector
\begin{equation}
\label{eq:aux_variables}
    \boldsymbol z^{(m)}(t)
    =
    (z^{(m)}_{1}(t), \cdots, z^{(m)}_{r_m}(t)),
\end{equation}
where $r_m$ is the number of forcing channels that represent the forcing term at order $m$.
The coupling vectors are columns of a matrix
\begin{equation}
    V^{(m)}
    =
    \begin{bmatrix}
    \boldsymbol v^{(m)}_{1}, \cdots, \boldsymbol v^{(m)}_{r_m}
    \end{bmatrix}
    \in\mathbb C^{J\times r_m},
\end{equation}
such that $\boldsymbol f^{(m-1)}(t)
    =\sum_{r=1}^{r_m} \boldsymbol v^{(m)}_{r}\, z^{(m)}_{r}(t)$.
With Eqs. \eqref{eq:single_channel_force_f} and \eqref{eq:zm_evolution}, a compact representation of the forcing term is obtained, where $t$ is the variable and $\Lambda^{(m)}=\operatorname{diag}(\lambda^{(m)}_{1},\ldots,\lambda^{(m)}_{r_m})$.

The solution vector $\boldsymbol u(\boldsymbol x, t)$ in Eq. \eqref{eq:mu_def_f} is now replaced by $\boldsymbol c(t)$.
The recursive system of deformation equations can then be written as one linear homogeneous autonomous system $d\boldsymbol y(t)/dt=A\,\boldsymbol y(t)$, which is defined with a lower block triangular matrix
\begin{equation}
    A=
    \begin{bmatrix}
        M & 0 & 0 & 0 & \cdots & 0 & 0\\
        0 & -\Lambda^{(1)} & 0 & 0 & \cdots & 0 & 0\\
        0 & V^{(1)} & M & 0 & \cdots & 0 & 0\\
        0 & 0 & 0 & -\Lambda^{(2)} & \cdots & 0 & 0\\
        0 & 0 & 0 & V^{(2)} & \cdots & 0 & 0\\
        \vdots & \vdots & \vdots & \vdots & \ddots & -\Lambda^{(\tilde{m})} & 0\\
        0 & 0 & 0 & 0 & \cdots & V^{(\tilde{m})} & M
    \end{bmatrix}.
    \label{eq:A_block_form}
\end{equation}
In Eq. \eqref{eq:A_block_form}, $M$ is the matrix representation of the linear operator. 
Each diagonal block $-\Lambda^{(m)}$ contains the decay rates of the exponential forcing channels introduced at order $m$. Each block $V^{(m)}$ maps the auxiliary forcing channels into the corresponding modal coefficients $\boldsymbol c^{(m)}$.

\subsection{Error bound for LHAM simulation}

In this section, the upper-bound of LHAM error is estimated. The total error for solving PDE is decomposed into discretization error, homotopy error, and dilation error. Lemma \ref{lem:ding_lindblad_simulation}, as follows, is adapted from Theorem $1$ of Ref. \cite{ding2024simulating} for Stinespring dilation error.

\begin{lemma}[Stinespring dilation error]
\label{lem:ding_lindblad_simulation}
Let $\mathcal L$ be a time-independent Lindbladian with system Hamiltonian $H$ and $K$ jump operators $\{F_k\}_{k=1}^{K}$ defined in Eq. \eqref{eq:lindbladian}.
The induced trace norm of the Lindbladian $\|\mathcal{L}\|_{1\to1}$ is bounded by an upper-bound parameter 
\begin{equation}
\label{eq:lindbladian_upper_bound}
    \mathscr B_{J,\tilde{m}}
    :=
    1+\|H\|_2+\sum_{k=1}^{K}\|F_k\|^2_2,
\end{equation}
since $\|\mathcal{L}\|_{1\to1} := \mathrm{sup}_{\|X\|_1 = 1}\,\|\mathcal{L}(X)\|_1 \le 2\|H\|_2+2\sum_{k=1}^{K}\|F_k\|_2^2  \le 2\mathscr B_{J,\tilde{m}}$.
For every integer $S\geq1$, there exist constants $c_{S,K},C_{S,K}>0$, depending only on $S$ and $K$, and a Hermitian dilated Hamiltonian in Eq.~\eqref{eq:ding_dilated_hamiltonian}, whenever $
    \Delta t\,\mathscr B_{J,\tilde m}\leq c_{S,K}$,
the $S$th-order approximation of the quantum channel
$\mathcal E_{\Delta t}^{(S)}$ in Eq. \eqref{eq:stinespring_dilation} defined with $\ket{0}_E = \ket{0^{a_S}}$ satisfies
\begin{equation}
\label{eq:ding_local_error}
    \sup_{\rho\in\mathbb D}
    \left\|
        e^{\mathcal L\Delta t}(\rho)
        -
        \mathcal E_{\Delta t}^{(S)}(\rho)
    \right\|_1
    \leq
    C_{S,K}
    \mathscr B_{J,\tilde m}^{S+1}
    \Delta t^{S+1},
\end{equation}
where $\mathbb D$ is the set of density operators.
If $\tilde n=T/\Delta t\in\mathbb N$, then the accumulated error at time $T$ satisfies
\begin{equation}
\label{eq:ding_global_error}
    \left\|
        e^{\mathcal LT}(\rho_0)
        -
        \left(\mathcal E_{\Delta t}^{(S)}\right)^{\tilde n}(\rho_0)
    \right\|_1
    \leq
    C_{S,K}T
    \mathscr B_{J,\tilde m}^{S+1}\Delta t^S.
\end{equation}
The general construction uses at most $
    a_S
    \leq
    \left\lceil
        (S+1)\log_2(K+1)
    \right\rceil$
environment ancilla qubits.
\end{lemma}

\begin{remark}
In the Lindblad operator embedding of LHAM, Eq. \eqref{eq:hamiltonian_jumpop} is equivalent to $H=\hat{A}_2\oplus0$ and $F=\sqrt{2\hat{A}_1}\oplus0$.
It corresponds to $K=1$ in Eq. \eqref{eq:lindbladian_upper_bound}. 
Since the zero blocks do not change the operator norm,
$
\mathscr B_{J,\tilde{m}}
=
1
+
\|\hat{A}_{2}\|_2
+
\|\sqrt{2\hat{A}_{1}}\|^2_2$.
As $\hat{A}_{1}\succeq 0$, this becomes $
    \mathscr B_{J,\tilde{m}}
    =
    1
    +
    \|\hat{A}_{2}\|_2
    +
    2\|\hat{A}_{1}\|_2$.
In addition, $C_S:=C_{S,1}>0$ is the dilation truncation constant for the order-$S$, single-jump operator dilation and is independent of $\Delta t$.
\end{remark}

\begin{remark}
For the explicit constructions in Eq.~\eqref{eq:ding_dilated_hamiltonian}, the environment Hilbert space is $
    \mathbb H_E^{(S)}
    =
    \operatorname{span}
    \{|0\rangle,\ldots,|S\rangle\}$, where $S\in\{1,2\}$,
and therefore has dimension $S+1$. Consequently, the required number
of environment qubits is $
    a_S
    =
    \left\lceil\log_2(S+1)\right\rceil$.
Hence, the first- and second-order constructions require one and two environment qubits, respectively.
\end{remark}

The Stinespring dilation error for LHAM is shown in Theorem \ref{thm:stinespring-dilation_error}, as follows. The proof of Theorem \ref{thm:stinespring-dilation_error} is provided in Appendix \ref{sec:proof2}.

\begin{theorem}[Stinespring dilation error for LHAM]
\label{thm:stinespring-dilation_error}
Let $\boldsymbol u_{J,\tilde m}(T)$ be the $\tilde m$-th order LHAM approximation after truncation with cutoff $J$, and let $\widehat{\boldsymbol u}_{J,\tilde m,S,\Delta t}(T)$ be the solution recovered from the upper-right block $\rho_{01}$ after applying the $S$-th order Stinespring dilation scheme with time step $\Delta t$. Let $\tilde n:=T/\Delta t\in\mathbb N$ be the number of time steps, and let $\mathcal E_{\Delta t}^{(S)}$ denote the one-step $S$-th order Stinespring dilation to approximate the Lindbladian channel
$e^{\mathcal L\Delta t}$.
Assume that the LHAM Lindbladian embedding contains one jump operator and that $
    \Delta t\,\mathscr B_{J,\tilde m}
    \leq
    c_{S,1}$, where $c_{S,1}>0$ is the small-step constant associated with the $S$-th order, single-jump operator dilation. Then the dilation error is
\begin{equation}
\label{eq:dilation_error}
\begin{aligned}
    \varepsilon_d
    &:=
    \left\|
        \boldsymbol u_{J,\tilde m}(T)
        -
        \widehat{\boldsymbol u}_{J,\tilde m,S,\Delta t}(T)
    \right\|_{L^2(\Omega)}
    \\
    &\leq
    2C_S\|\boldsymbol y(0)\|_2
    \sqrt{\tilde m+1}\,
    e^{\gamma T}
    T\mathscr B_{J,\tilde m}^{S+1}
    \Delta t^S.
\end{aligned}
\end{equation}
\end{theorem}

\begin{remark}
The preceding estimate accounts for the truncation error of the order-$S$ Stinespring dilation. It assumes that the dilated unitary $U=e^{-i\sqrt{\Delta t}\widetilde{H}^{(S)}}$ is implemented exactly and that the recovered density-matrix block is known exactly. If the unitary is approximated by a Hamiltonian-simulation or product-formula algorithm, or if the block is estimated using finitely many measurements, the corresponding implementation and sampling errors must be added separately.
\end{remark}

Here, the Euclidean norm of the Fourier-coefficient error is identified with the spatial $L^2(\Omega)$-norm using an orthonormal Fourier basis and Parseval's identity.
The total error for solving PDE with LHAM is shown in Theorem \ref{thm:homotopy_dilation_balance}, as follows. The proof of Theorem \ref{thm:homotopy_dilation_balance} is provided in Appendix \ref{sec:proof3}.

\begin{theorem}[Total error for LHAM]
\label{thm:homotopy_dilation_balance}
Let $\boldsymbol u(T)$ be the exact solution of the target nonlinear evolution equation, and let $
    \boldsymbol u_J(T):=P_J\boldsymbol u(T)$
be its Fourier projection with cutoff $J$. Assume that
$
    \boldsymbol u
    \in
    \mathscr C\!\left([0,T];H^f(\Omega)\right)$, $f>0$,
so that
\begin{equation}
    \left\|
        \boldsymbol u(T)-\boldsymbol u_J(T)
    \right\|_{L^2(\Omega)}
    \leq
    C_fJ^{-f},
\end{equation}
where $C_f>0$ is independent of $J$.
Suppose that the homotopy-series solution $
    \boldsymbol u_J(T)
    =
    \sum_{m=0}^{\infty}
    \boldsymbol u_J^{(m)}(T)$
satisfies the geometric-decay condition
\begin{equation}
\label{eq:homotopy_geometric_decay}
    \left\|
        \boldsymbol u_J^{(m+1)}(T)
    \right\|_{L^2(\Omega)}
    \leq
    \theta_J
    \left\|
        \boldsymbol u_J^{(m)}(T)
    \right\|_{L^2(\Omega)},
    \qquad
    m\geq0,
\end{equation}
for some contraction factor $\theta_J\in(0,1)$. Define
$
    \boldsymbol u_{J,\tilde m}(T)
    :=
    \sum_{m=0}^{\tilde m}
    \boldsymbol u_J^{(m)}(T)$.
Let the lifted LHAM matrix $A$ in Eq.~\eqref{eq:A_block_form} be embedded into shifted Lindblad dynamics with shift $\gamma\geq0$, system Hamiltonian $H$, and one jump operator $F$ as in Eq.~\eqref{eq:hamiltonian_jumpop}. Let $\mathscr B_{J,\tilde m}$ be defined in Eq.~\eqref{eq:lindbladian_upper_bound}, and assume that $
    \Delta t\,\mathscr B_{J,\tilde m}
    \leq
    c_{S,1}$.
Then
\begin{equation}
\label{eq:lham_homotopy_dilation_error}
\begin{aligned}
\left\|
    \boldsymbol u(T)
    -
    \widehat{\boldsymbol u}_{J,\tilde m,S,\Delta t}(T)
\right\|_{L^2(\Omega)}
\leq\;&
C_fJ^{-f}
+
\left\|
    \boldsymbol u_J^{(0)}(T)
\right\|_{L^2(\Omega)}
\frac{\theta_J^{\tilde m+1}}{1-\theta_J}
\\
&+
2C_S\|\boldsymbol y(0)\|_2
\sqrt{\tilde m+1}\,
e^{\gamma T}
T\mathscr B_{J,\tilde m}^{S+1}
\Delta t^S.
\end{aligned}
\end{equation}
\end{theorem}

The choices of Fourier cutoff, homotopy order, dilation order, and time step are guided as follows.
\begin{corollary}
\label{cor:choice_M_q_dt}
Let $\varepsilon>0$ be a target error tolerance. It is sufficient to allocate $\varepsilon/3$ to each of the three error terms:
\begin{align}
    C_fJ^{-f}
    &\leq
    \frac{\varepsilon}{3},
    \label{eq:target_fourier_error}
    \\
    \left\|
        \boldsymbol u_J^{(0)}(T)
    \right\|_{L^2(\Omega)}
    \frac{\theta_J^{\tilde m+1}}{1-\theta_J}
    &\leq
    \frac{\varepsilon}{3},
    \label{eq:target_homotopy_error}
    \\
    2C_S\|\boldsymbol y(0)\|_2
    \sqrt{\tilde m+1}\,
    e^{\gamma T}
    T\mathscr B_{J,\tilde m}^{S+1}
    \Delta t^S
    &\leq
    \frac{\varepsilon}{3}.
    \label{eq:target_dilation_error}
\end{align}
A sufficient Fourier cutoff is
\begin{equation}
\label{eq:sufficient_fourier_cutoff}
    J
    \geq
    \left\lceil
        \left(
            \frac{3C_f}{\varepsilon}
        \right)^{1/f}
    \right\rceil.
\end{equation}
Assuming $\|\boldsymbol u_J^{(0)}(T)\|_{L^2(\Omega)}>0$, a sufficient homotopy order is
\begin{equation}
\label{eq:sufficient_homotopy_order}
    \tilde m
    \geq
    \max\left\{
        0,\,
        \left\lceil
        \frac{
            \log\!\left[
                \dfrac{
                    \varepsilon(1-\theta_J)
                }{
                    3\|
                    \boldsymbol u_J^{(0)}(T)
                    \|_{L^2(\Omega)}
                }
            \right]
        }{
            \log\theta_J
        }
        -1
        \right\rceil
    \right\}.
\end{equation}
For a fixed time step satisfying $
    0<
    \Delta t\,\mathscr B_{J,\tilde m}
    <1$,
suppose that $C_r\leq\overline C$ for every candidate dilation order $r$ under consideration, where $\overline C>0$ is independent of $r$. A sufficient dilation order is
\begin{equation}
\label{eq:sufficient_dilation_order}
    S
    \geq
    \max\left\{
        1,\,
        \left\lceil
        \frac{
            \log\!\left[
                \dfrac{\varepsilon}{
                    6\overline C
                    \|\boldsymbol y(0)\|_2
                    \sqrt{\tilde m+1}\,
                    e^{\gamma T}
                    T\mathscr B_{J,\tilde m}
                }
            \right]
        }{
            \log\!\left(
                \Delta t\,
                \mathscr B_{J,\tilde m}
            \right)
        }
        \right\rceil
    \right\}.
\end{equation}
Alternatively, for a fixed dilation order $S$, it is sufficient to choose
\begin{equation}
\label{eq:sufficient_dilation_timestep}
    \Delta t
    \leq
    \min\left\{
        \frac{c_{S,1}}{\mathscr B_{J,\tilde m}},
        \left[
        \frac{\varepsilon}{
            6C_S
            \|\boldsymbol y(0)\|_2
            \sqrt{\tilde m+1}\,
            e^{\gamma T}
            T\mathscr B_{J,\tilde m}^{S+1}
        }
        \right]^{1/S}
    \right\},
\end{equation}
with $\tilde n=T/\Delta t\in\mathbb N$.
\end{corollary}

\begin{remark}
The first term in Eq.~\eqref{eq:lham_homotopy_dilation_error} represents the spatial discretization error. For a general spatial discretization, let $\boldsymbol u_h$ denote the spatially discrete solution, where $h>0$ is a spatial-resolution parameter. Suppose that
\begin{equation}
    \left\|
        \boldsymbol u(T)-\boldsymbol u_h(T)
    \right\|_{L^2(\Omega)}
    \leq
    C_s\Gamma_s(h),
\end{equation}
where $\Omega$ is the spatial domain and $C_s>0$ is independent of $h$ and $\Gamma_s(h)\to0$ as $h\to0$. For the Fourier projection considered above, $h=\Theta(J^{-1})$ and $\Gamma_s(h)=h^f$, equivalently, $\Gamma_s=\Theta(J^{-f})$
\cite{shen2011spectral,Hesthaven_Gottlieb_Gottlieb_2007}.
For a stable finite-difference discretization known to converge with spatial order $p$, one has $\Gamma_s(h)=h^p$, provided that the exact solution has the required regularity \cite{leveque2007finite}.
An analogous total-error estimate follows by replacing the Fourier truncation term and all other spatial-discretization-dependent quantities with those corresponding to the chosen discretization.
\end{remark}

Consider an $s$-component field on a $d_x$-dimensional periodic domain, where $s$ is the number of dependent variables. A tensor-product Fourier discretization with cutoff $|j_k|\leq J$,
$k=1,\ldots,d_x$, has $
    D_J=s(2J+1)^{d_x}
       =\Theta\!\left(sJ^{d_x}\right)$
spatial degrees of freedom. If each of $\mathscr O(D_J)$ spatial channels requires temporal factors $
    t^r e^{-\lambda t}$,
where $r=0,\ldots,\tilde m$, and $\lambda$ is a decay rate and $r$ is the temporal monomial degree, then the lifted dimension satisfies $
    \mathscr O\!\left(D_J(\tilde m+1)\right)$.

For comparison, consider a uniform finite-difference mesh with spacing $h$ on a domain of fixed size. If $N_h=\Theta(h^{-1})$ is the number of grid points per spatial direction, then its number of spatial degrees of freedom satisfies $
    D_h
    =
    \Theta\!\left(sN_h^{d_x}\right)
    =
    \Theta\!\left(sh^{-d_x}\right)$.
Therefore, Fourier and finite-difference discretizations have the same asymptotic dimension scaling with the number of degrees of freedom per spatial direction, although their exact dimensions and their convergence rates with respect to spatial resolution generally differ. Under the same temporal-channel construction, the
finite-difference-based lifted dimension similarly scales as $\mathscr O(D_h(\tilde m+1))$.

Increasing the spatial resolution reduces the spatial discretization error but enlarges the spatial and lifted systems. Increasing the homotopy order $\tilde m$ reduces the homotopy truncation error but also enlarges the lifted block matrix $A$. Since the Hamiltonian $H$, the jump operators $F_k$, the shift parameter $\gamma$, and $\mathscr B_{J,\tilde m}$ are constructed from $A$, their norms or bounds may increase with $\tilde m$, although monotonic growth is not guaranteed. Consequently, the dilation-error bound and simulation cost may increase even while the homotopy truncation error decreases.
Thus, the appropriate value of $\tilde m$ need not be the largest feasible order. It should balance the homotopy truncation error against $\mathscr B_{J,\tilde m}$, $\gamma$, the dilation error, and the simulation cost.

\section{\label{sec:demonstration}Demonstrations}
\subsection{Burgers' Equation}

Burgers' equation describes the nonlinear wave propagation and shock formation.
The one-dimensional viscous Burgers' equation is defined as 
\begin{equation}
    \frac{\partial}{\partial t} u(x,t)
    = \nu\,\frac{\partial^2}{\partial x^2}u(x,t) - u(x,t)\,\frac{\partial}{\partial x} u(x,t),
    \label{eq:burgers_pde}
\end{equation}
including the position variable $x\in[0,2\pi)$, time $\ t\ge 0$, and the kinematic viscosity $\nu>0$. $u(x,t)$ is the velocity field with periodic boundary conditions. The initial condition is set as $u(x,0)=u_0(x)$. 
The linear and nonlinear components in Eq. \eqref{eq:burgers_pde} are the diffusion and advection terms, respectively.
The Hilbert space is defined as $\mathbb{H} = L^2(\Omega)$ on the periodic spatial domain $\Omega=[0,2\pi)$.

After the functional expansion in the Fourier basis with $|j| \le 4$, the linear homogeneous autonomous system is defined with the block matrix $A$ in Eq. \eqref{eq:A_block_form}. Details of the derivation is provided in Appendix \ref{sec:Burgers}.
The initial condition was chosen as $u(x,0)=\sin x$.
The homotopy-series solution was truncated at $\tilde{m}=4$.
The results from LHAM and finite difference method (FDM) were compared with each other.
The difference between the solutions $u^\mathrm{L}$ from LHAM  and $u^\mathrm{C}$ from FDM is quantified with
relative $L^2$ norm error
$\varepsilon_{L^2}
=
\sqrt{
\sum_{i=1}^{N}
\left(
u_i^{\mathrm{L}}
-
u_i^{\mathrm{C}}
\right)^2
/
\sum_{i=1}^{N}\left(
u_i^{\mathrm{C}}
\right)^2 } 
$, where $N$ is the number of nodes.

\begin{figure}[h]
    \centering
    \includegraphics[width=0.45\textwidth]{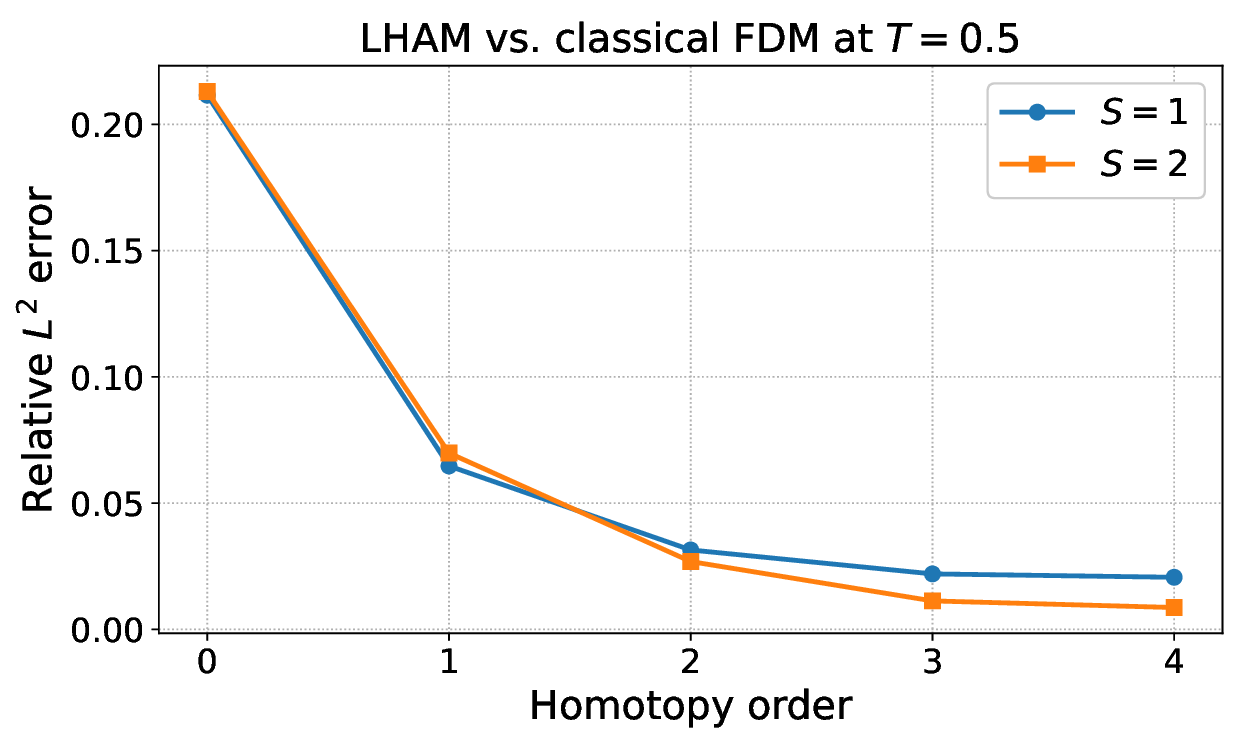}
    \caption{Relative $L^2$ norm error of solving Burgers' equation with LHAM for the first- and second-order Stinespring dilation}
    \label{fig:burgers_L2}
\end{figure}

\begin{figure}[h]
    \centering
    \includegraphics[width=0.45\textwidth]{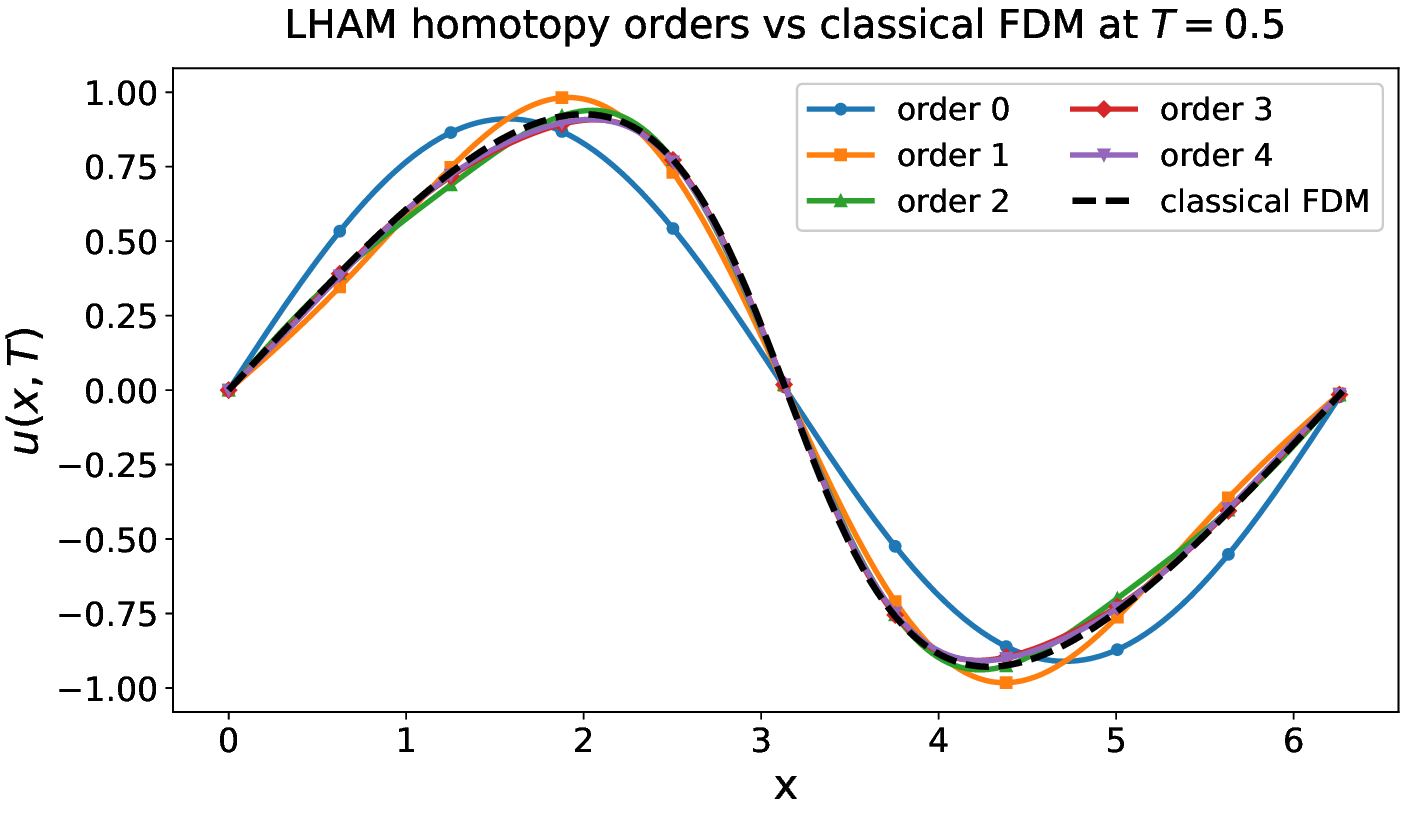}
    \caption{Calculated $u$ values from Burgers' equation with FDM and the first-order dilation LHAM at different truncation orders}
    \label{fig:burgers_profiles_o1}
\end{figure}

\begin{figure}[h]
    \centering
    \includegraphics[width=0.45\textwidth]{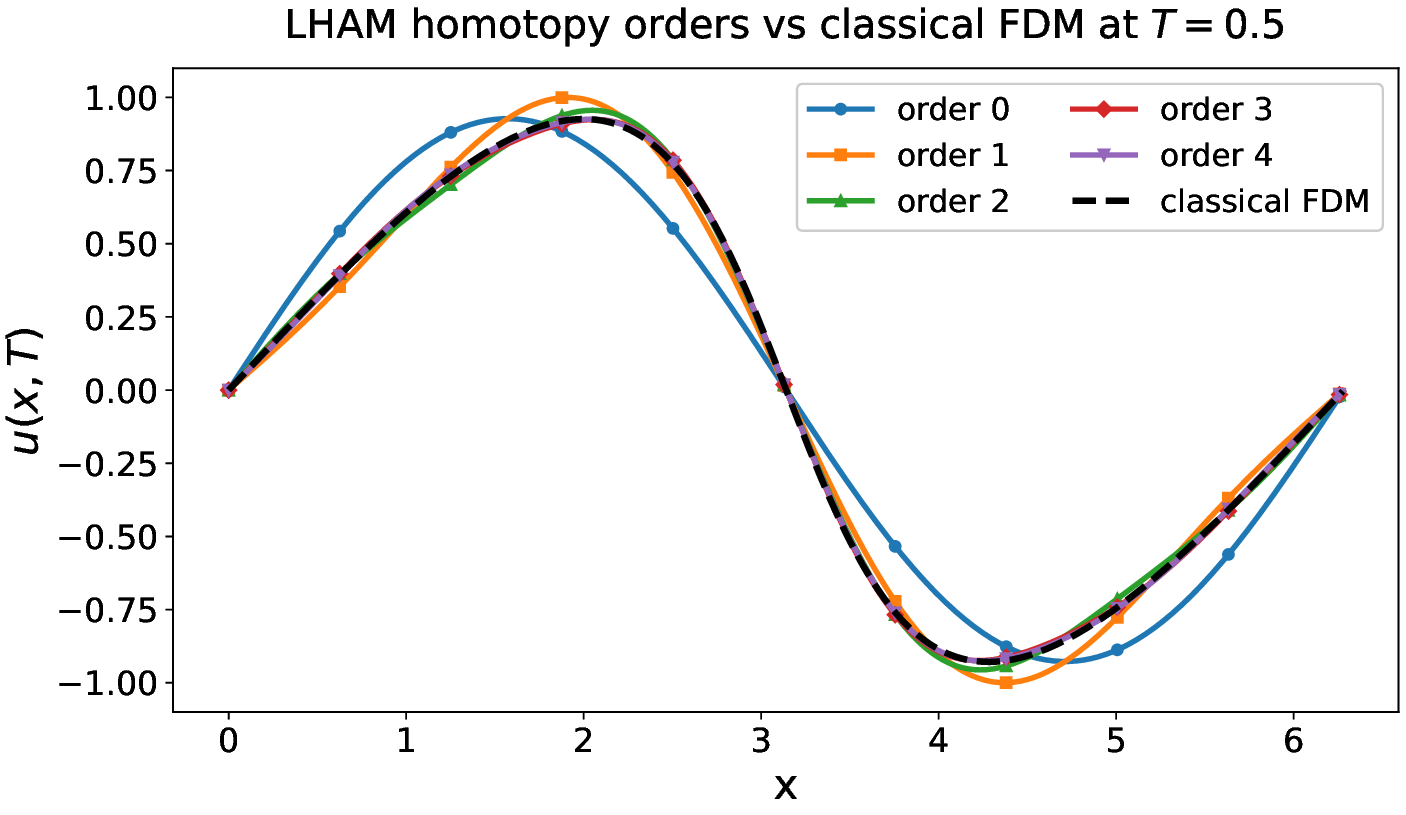}
    \caption{Calculated $u$ values from Burgers' equation with FDM and the second-order dilation LHAM at different truncation orders}
    \label{fig:burgers_profiles_o2}
\end{figure}

\begin{figure*}[t]
    \centering
    \includegraphics[width=0.5\textwidth]{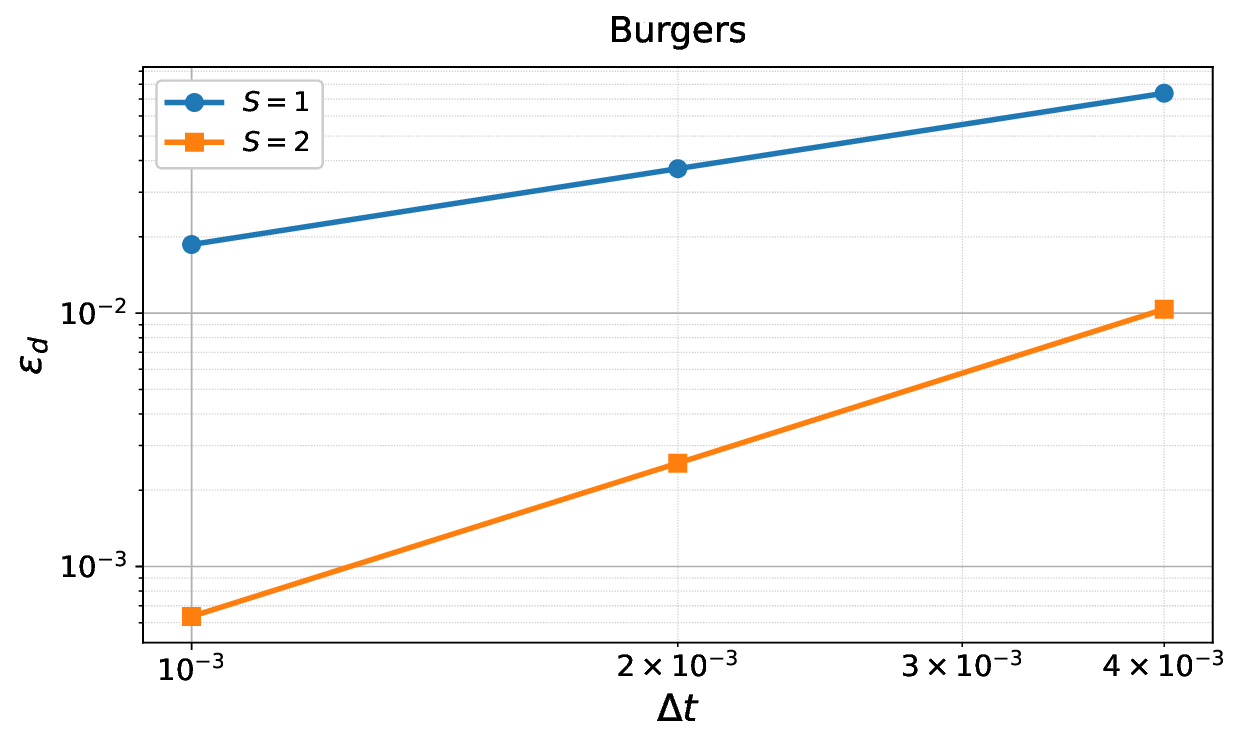}
    \caption{Dilation error of solving Burgers' equation with LHAM for the first- and second-order Stinespring dilation. The theoricial slope for $\mathscr O(\Delta t^1)$ and $\mathscr O(\Delta t^2)$ lines are 1 and 2, respectively, in the logarithmic scale.}
    \label{fig:Burgers_dilation_error}
\end{figure*}

Each time step was set to $\Delta t=0.001$ and the quantum state which encodes the solution is evolved for $500$ time steps for a total evolution time $T=0.5$.
The kinematic viscosity was set to $\nu=0.15$. 
The resulting errors are plotted in Fig.  \ref{fig:burgers_L2}, where the LHAM converges as the homotopy order increases.
From the second homotopy order, the second-order dilation $S=2$ resulted in less error than the first-order dilation $S=1$. 
At the fourth homotopy order, the relative $L^2$ norm errors are $2.062\%$ and $0.864\%$ for $S=1$ and $S=2$, respectively.
The values of $u(x,T=0.5)$ calculated from LHAM and FDM are shown in Figs. \ref{fig:burgers_profiles_o1} and \ref{fig:burgers_profiles_o2} for $S=1$ and $S=2$, respectively. It is observed that the solution from LHAM closely approximates the solution from FDM starting from $\tilde{m}=3$ for $S=2$.
The Stinespring dilation error $\varepsilon_d$ is calculated by varing time steps $\Delta t$ as $0.001$, $0.002$, and $0.004$ in Fig. \ref{fig:Burgers_dilation_error}. Theoretically, the slopes for $S=1$ and $S=2$ are $1$ and $2$ in the logarithmic scale as the dilation error scales as $\mathscr O(\Delta t^1)$ and $\mathscr O(\Delta t^2)$ based on Theorem \ref{thm:stinespring-dilation_error}, respectively. The measured slopes are $0.992$ and $2.014$ for $S=1$ and $S=2$, which match the predictions from Theorem \ref{thm:stinespring-dilation_error}.

\subsection{Magnetohydrodynamics}

Magnetohydrodynamics (MHD) is a model of electrically conducting fluids that treats charged particles as a continuous fluid.
The reduced MHD equations, which are applicable for incompressible fluids, are described by a coupled system of nonlinear PDEs defined as
\begin{equation}
\begin{aligned}
\frac{\partial \omega}{\partial t}
&=
\nu \nabla^2 \omega
-\left(
\frac{\partial \phi}{\partial x}\frac{\partial \omega}{\partial y}
-
\frac{\partial \phi}{\partial y}\frac{\partial \omega}{\partial x}
\right)
+\left(
\frac{\partial \xi}{\partial x}\frac{\partial \zeta}{\partial y}
-
\frac{\partial \xi}{\partial y}\frac{\partial \zeta}{\partial x}
\right), \\
\frac{\partial \xi}{\partial t}
&=
\eta \nabla^2 \xi
-\left(
\frac{\partial \phi}{\partial x}\frac{\partial \xi}{\partial y}
-
\frac{\partial \phi}{\partial y}\frac{\partial \xi}{\partial x}
\right),
\end{aligned}
\label{eq:rmhd_pde}
\end{equation}
where $(x,y)\in[0,2\pi)^2$ are spatial coordinates with periodic boundary conditions, $\nu>0$ is the kinematic viscosity, and $\eta>0$ is the resistivity. 
The vorticity field $\omega(x,y,t)$, the magnetic potential $\xi(x,y,t)$, and the stream function $\phi(x,y,t)$ are related through
\begin{equation}
\omega=-\nabla^2\phi,
\qquad
\zeta=-\nabla^2\xi,
\end{equation}
where $\zeta(x,y,t)$ is the current density. The linear components of Eq.~\eqref{eq:rmhd_pde} are the diffusive terms $\nu\nabla^2\omega$ and $\eta\nabla^2\xi$, while the nonlinear components describe fluid advection and magnetic coupling.
The Hilbert space is defined as $\mathbb{H}=L^2(\Omega)\times L^2(\Omega)$ on the periodic two-dimensional domain $\Omega=[0,2\pi)^2$.

After the functional expansion in Fourier basis with $|j| \le 4$, the linear homogeneous autonomous system is defined with the block matrix $A$ in Eq. \eqref{eq:A_block_form}. Details of the derivation are provided in Appendix \ref{sec:MHD}.
The initial conditions are given by
\begin{align}
\label{eq:MHD_IC}
\begin{split}
\omega(x,y,0) &= \sin x+\frac12\sin(x-y), \\
\xi(x,y,0) &= \cos y+\frac14\cos(x+y).
\end{split}
\end{align}
The homotopy-series solution was truncated at $\tilde{m}=4$. 
The kinematic viscosity and resistivity were set to $\nu=0.3$ and $\eta=0.3$, respectively.
Each time step was set $\Delta t =0.005$ and the state was evolved over $100$ time steps for a total evolution time of $T=0.5$.
The combined relative $L^2$ norm error was quantified as
$
\varepsilon_{L^2}^{(\mathrm{comb})}
=
\sqrt{
\sum_{i=1}^{N}
\bigl(\left(
\omega_i^{\mathrm{L}} - \omega_i^{\mathrm{C}}
\right)^2
+
\left(
\xi_i^{\mathrm{L}} - \xi_i^{\mathrm{C}}
\right)^2\bigr)
/
\sum_{i=1}^{N}
\bigl(\left(
\omega_i^{\mathrm{C}}
\right)^2
+
\left(
\xi_i^{\mathrm{C}}
\right)^2\bigr)
}
$.

\begin{figure}[h]
    \centering
    \includegraphics[width=0.45\textwidth]{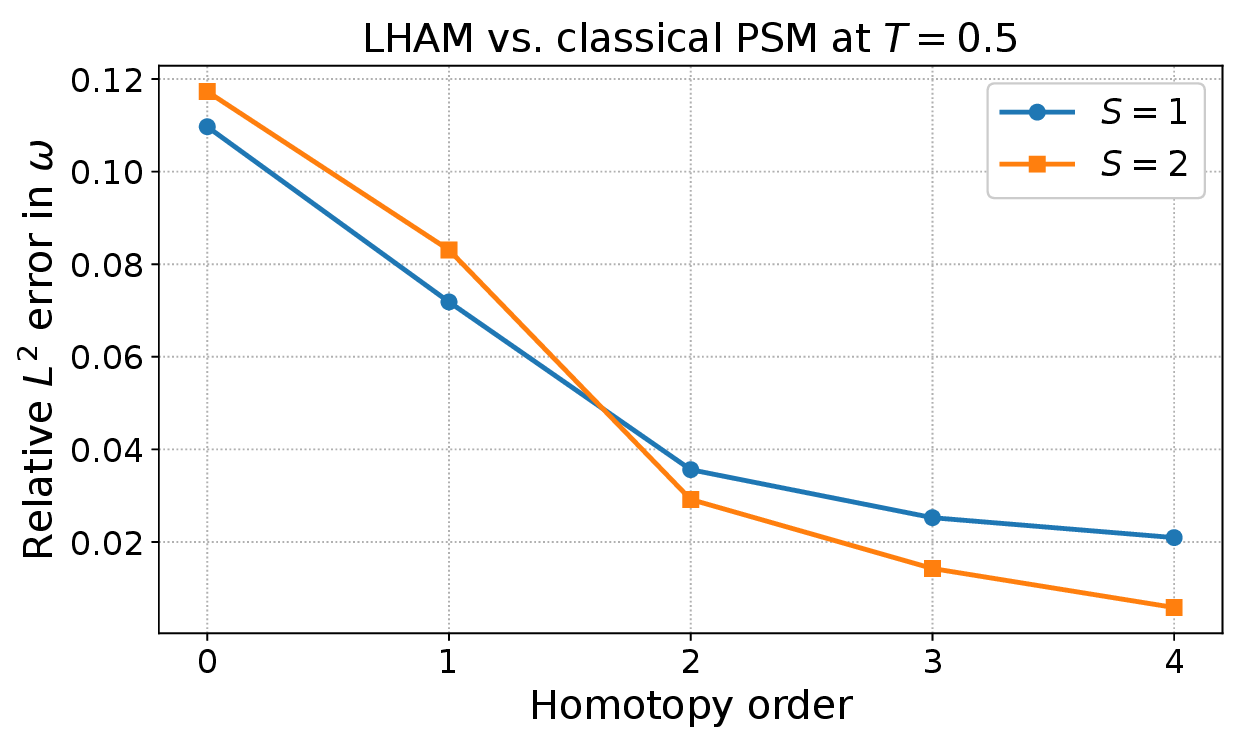}
    \caption{$\omega$ relative $L^2$ norm error of solving reduced MHD equations with LHAM for the first- and second-order Stinespring dilation.}
    \label{fig:MHD_L2_omega}
\end{figure}

\begin{figure}[h]
    \centering
    \includegraphics[width=0.45\textwidth]{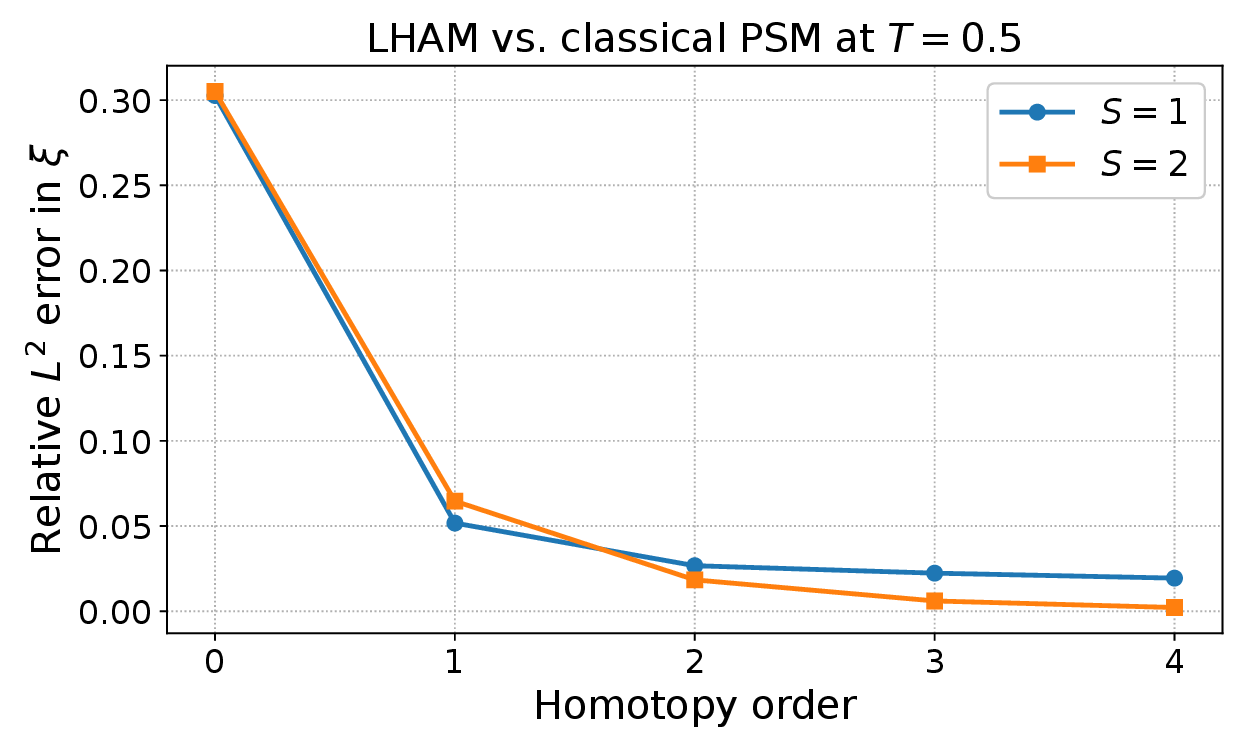}
    \caption{$\xi$ relative $L^2$ norm error of solving reduced MHD equations with LHAM for the first- and second-order Stinespring dilation.}
    \label{fig:MHD_L2_xi}
\end{figure}

\begin{figure}[h]
    \centering
    \includegraphics[width=0.45\textwidth]{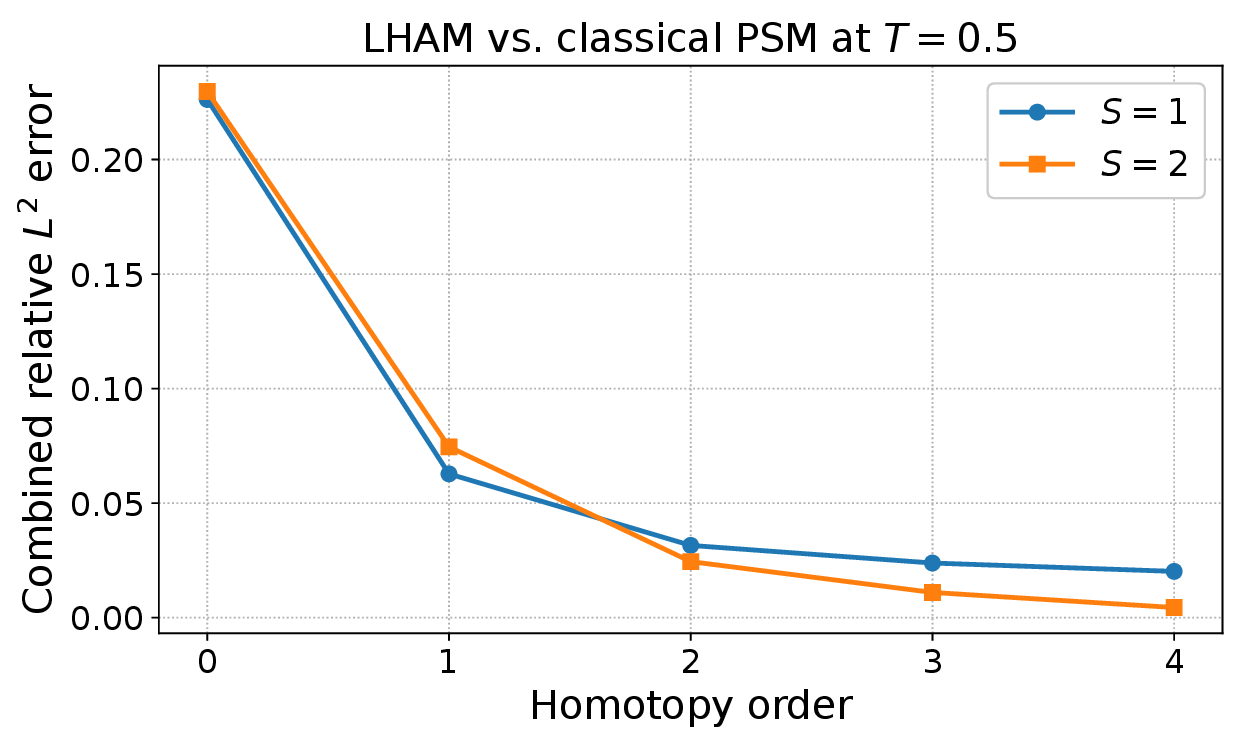}
    \caption{Combined relative $L^2$ norm error of solving reduced MHD equations with LHAM for the first- and second-order Stinespring dilation.}
    \label{fig:MHD_L2_combined}
\end{figure}

The results of LHAM and pseudo-spectral method (PSM) are compared with each other.
The relative $L^2$ norm errors are shown in Figs. \ref{fig:MHD_L2_omega}, \ref{fig:MHD_L2_xi}, and \ref{fig:MHD_L2_combined} for $\omega$, $\xi$, and the combined field, respectively. 
The second-order dilation errors showed less values than the first-order errors from the second homotopy order.
At $\tilde{m}=4$, the relative $L^2$ errors for $\omega$ are $2.093\%$ and $0.583\%$ for $S=1$ and $S=2$, respectively.
At $\tilde{m}=4$, the $\xi$ errors are $1.941\%$ and $0.215\%$ for $S=1$ and $S=2$.
The combined relative $L^2$ norm errors at $\tilde{m}=4$ were $2.020\%$ and $0.442\%$ for $S=1$ and $S=2$.
Overall, the $L^2$ error of $\xi$ decreases faster than that of $\omega$.
This is because the first equation in Eq. \eqref{eq:rmhd_pde} for $\omega$ has more complex nonlinear coupling terms consisting of $\phi$, $\xi$, and $\zeta$. The second equation in Eq. \eqref{eq:rmhd_pde} has simpler coupling terms consisting of $\xi$ and $\phi$ which results in smoother solutions. 
The error of $\omega$ can be further reduced by increasing number of Fourier basis functions and homotopy order.

\begin{figure*}[t]
    \centering
    \includegraphics[width=0.9\textwidth]{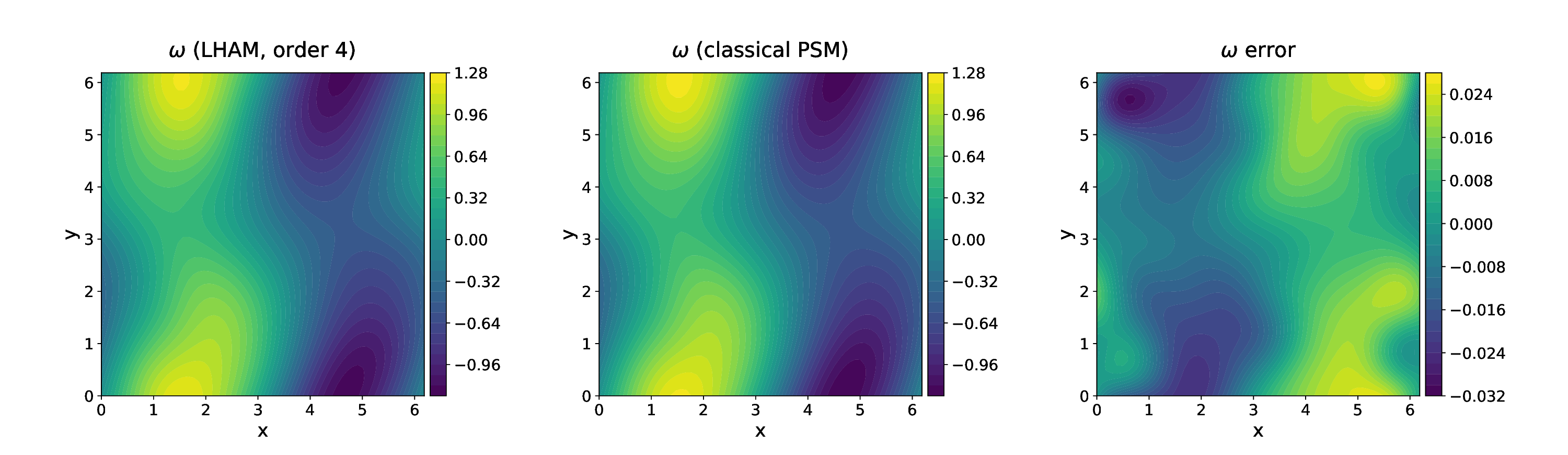}
    \caption{$\omega$ field profile of the first-order dilation LHAM and classical pseudo-spectral method results with the error.}
    \label{fig:MHD_omega_o1}
\end{figure*}

\begin{figure*}[t]
    \centering
    \includegraphics[width=0.9\textwidth]{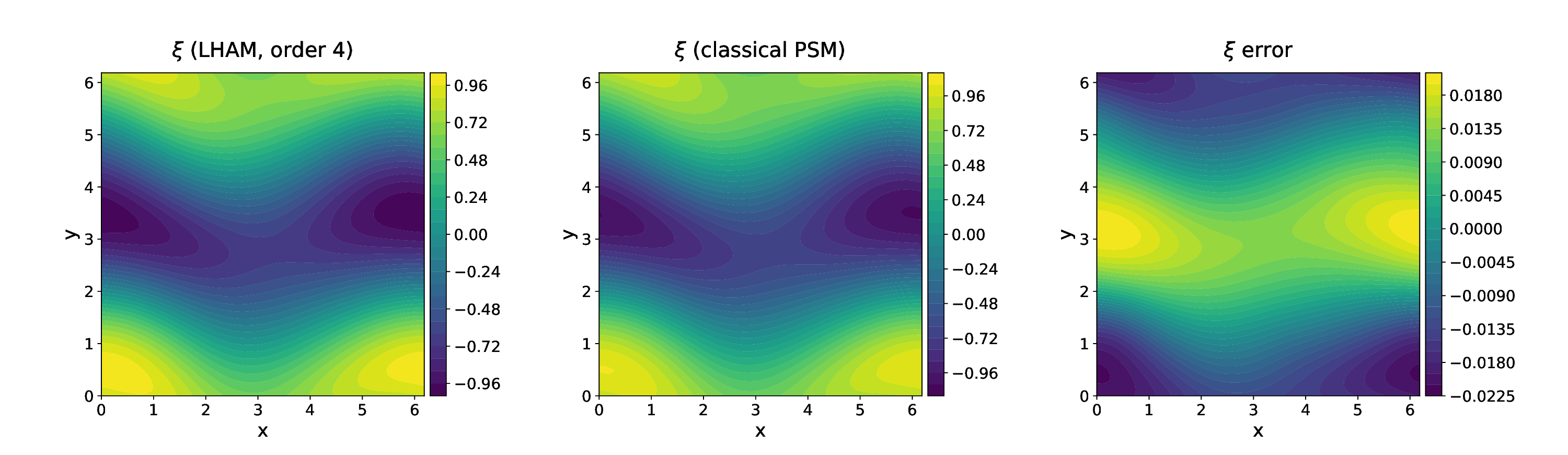}
    \caption{$\xi$ field profile of the first-order dilation LHAM and classical pseudo-spectral method results with the error.}
    \label{fig:MHD_psi_o1}
\end{figure*}

\begin{figure*}[t]
    \centering
    \includegraphics[width=0.9\textwidth]{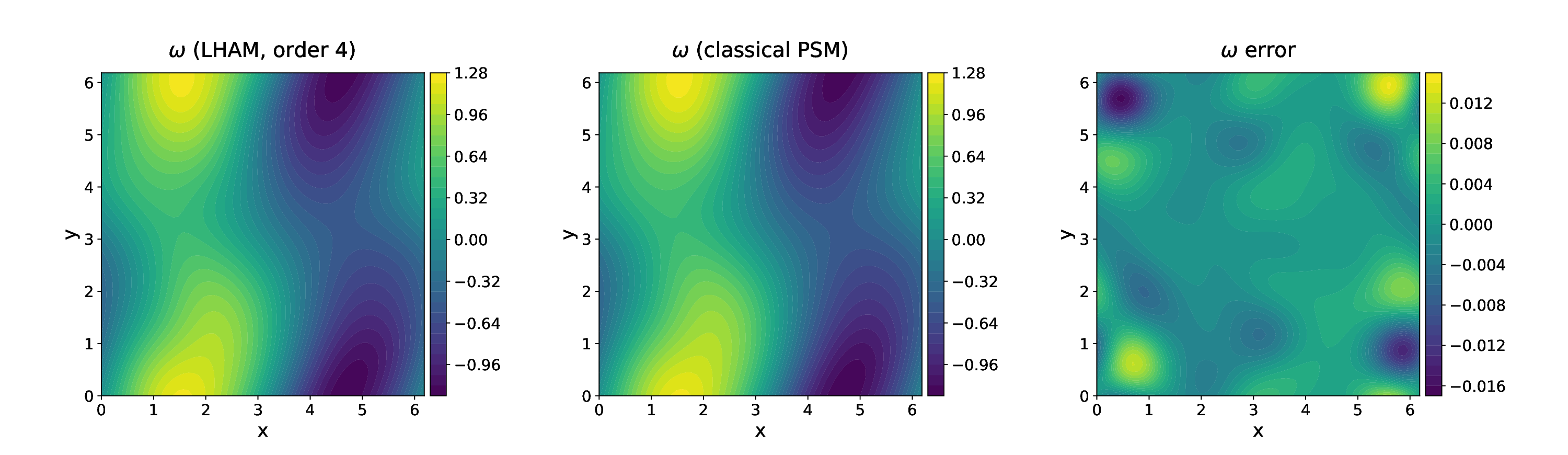}
    \caption{$\omega$ field profile of the second-order dilation LHAM and classical pseudo-spectral method results with the error.}
    \label{fig:MHD_omega_o2}
\end{figure*}

\begin{figure*}[t]
    \centering
    \includegraphics[width=0.9\textwidth]{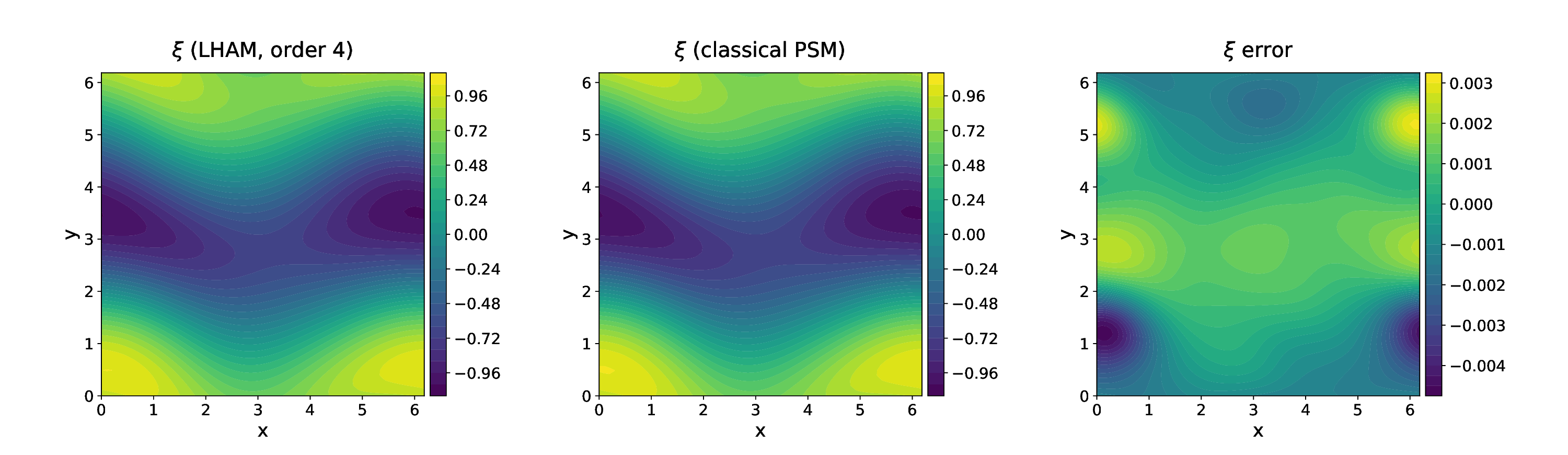}
    \caption{$\xi$ field profile of the second-order dilation LHAM and classical pseudo-spectral method results with the error.}
    \label{fig:MHD_psi_o2}
\end{figure*}

\begin{figure*}[t]
    \centering
    \includegraphics[width=0.5\textwidth]{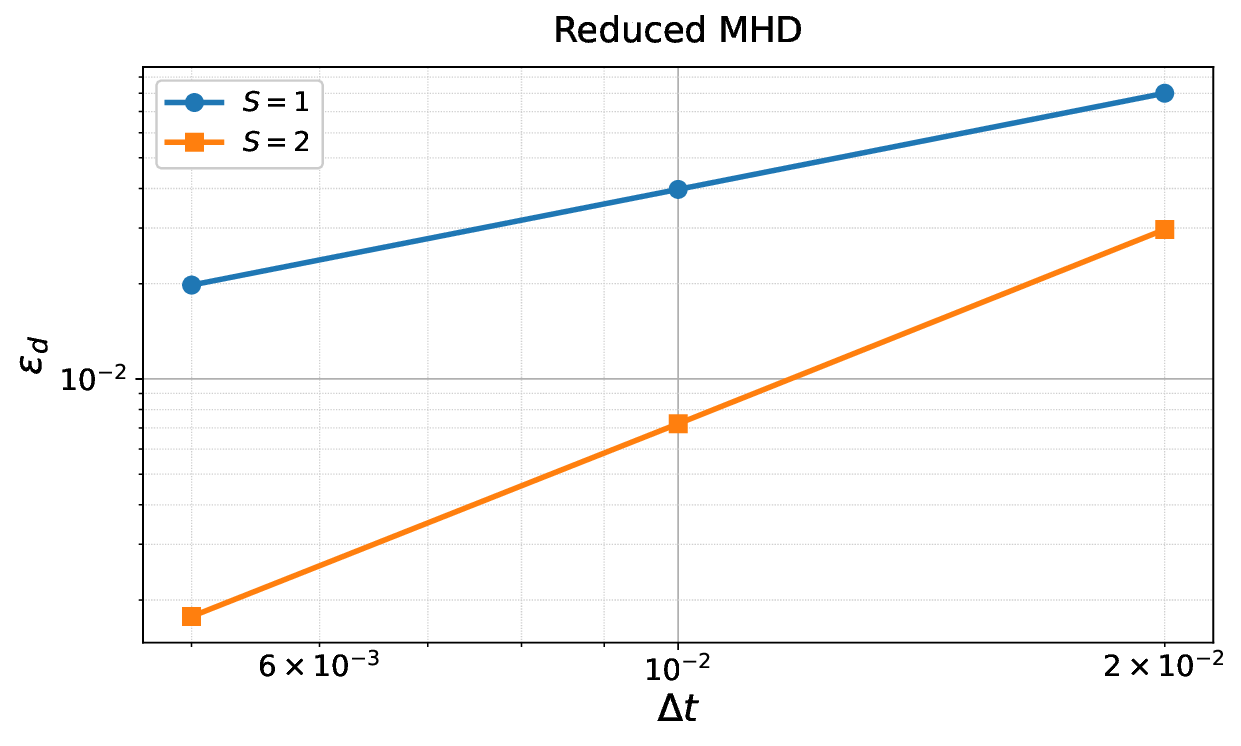}
    \caption{Dilation error of solving reduced MHD equations with LHAM for the first- and second-order Stinespring dilation. The theoricial slope for $\mathscr O(\Delta t^1)$ and $\mathscr O(\Delta t^2)$ lines are 1 and 2, respectively, in the logarithmic scale.}
    \label{fig:MHD_dilation_error}
\end{figure*}

The calculated values of $\omega(x,y,T=0.5)$ and $\xi(x,y,T=0.5)$ are shown in Figs. \ref{fig:MHD_omega_o1} and \ref{fig:MHD_psi_o1} for $S=1$ and Figs. \ref{fig:MHD_omega_o2} and \ref{fig:MHD_psi_o2} for $S=2$, respectively. It is observed that the solution from LHAM approximates the solution from PSM. The errors in Figs. \ref{fig:MHD_omega_o1}, \ref{fig:MHD_psi_o1}, \ref{fig:MHD_omega_o2}, and \ref{fig:MHD_psi_o2} were computed as $\varepsilon_D = u^\mathrm{L} - u^\mathrm{C}$.
The $\omega$ and $\xi$ errors are reduced by increasing the dilation order from $S=1$ to $S=2$ as well as the coherent spatial error patterns observed for both fields.
The oscillatory error patterns occur when higher Fourier modes are missing or have incorrect amplitudes.
The error plots indicate which Fourier modes are missing or have mismatched coefficients in LHAM compared to PSM.
The error can be further ruduced by increasing the number of Fourier modes and the dilation order.
The Stinespring dilation error $\varepsilon_d$ is calculated by varing time steps $\Delta t$ as $0.005$, $0.01$, and $0.02$ in Fig. \ref{fig:MHD_dilation_error}. The theoretical slopes $\mathscr O(\Delta t^1)$ and $\mathscr O(\Delta t^2)$ are plotted for $S=1$ and $S=2$, respectively. The measured slopes are $1.008$ and $2.032$ for $S=1$ and $S=2$, which are very close to the theoretical predictions.

\section{\label{sec:conclusion}Conclusions}

Most of the existing linearization-based quantum nonlinear ODE/PDE solvers face the scalability challenge since the state space dimension increases exponentially. In this paper, a new quantum method, LHAM, is proposed to solve nonlinear ODEs/PDEs. In LHAM, a linear homogeneous autonomous system is constructed and the dimension of the state space grows linearly with the homotopy order and spatial discretization.
The dimension of the Hilbert space is also exponentially smaller than the existing homotopy-based quantum methods in which the introduced auxiliary spatial variables or solution fields increase the state dimension combinatorially, because of the secondary linearization. 
To simulate nonunitary dynamics in the linearized PDEs, as few as two ancilla qubits are needed for the Lindblad dynamics simulation. 
In contrast, the required number of ancilla qubits increases logarithmically with the discretization of the auxiliary continuous variables in other quantum nonunitary dynamics solvers such as Sch\"odingerization or linear combination of Hamiltonian simulation.
Therefore, LHAM can significantly improve the scalability of quantum nonlinear differential equation solvers. 
The dissipative dynamics of the Lindblad master equation embedded in the jump operators can be generalized to encode different types of ODEs and PDEs.
The solution vector is recovered through measurements of Hermitian observables, and this process does not rely on post-selection success probability.

The proposed LHAM was used to solve Burgers' equation and the coupled two-dimensional MHD equations.
In Burgers' equation, the relative $L^2$ norm errors were $2.062\%$ and $0.864\%$ for $S=1$ and $S=2$, respectively, at the fourth homotopy order. 
The measured slopes for the Stinespring dilation error are $0.992$ and $2.014$ for $S=1$ and $S=2$, which are nearly identical with the theoretical slopes for $\mathscr O(\Delta t^1)$ and $\mathscr O(\Delta t^2)$.
In MHD, the vorticity and magnetic potential fields were closely approximated with the PSM solutions. The combined relative $L^2$ norm errors at $\tilde{m}=4$ were $2.020\%$ and $0.442\%$ for $S=1$ and $S=2$, respectively. 
The measured dilation error slopes are $1.008$ and $2.032$ for $S=1$ and $S=2$, which are nearly identical to the theoretical slopes for $\mathscr O(\Delta t^1)$ and $\mathscr O(\Delta t^2)$.
The solution accuracy can be further improved by increasing the homotopy order and spatial discretization. 
In addition, the measured Stinespring dilation error are matched with the theoretical prediction for LHAM.
The numerical accuracy of the mid-circuit-measure-and-reset scheme can also be improved by introducing the higher-order term in the Taylor expansion of the nonunitary time-evolution operator \cite{ding2024simulating}.

The LHAM was demonstrated with second-order nonlinear PDEs. Future work will include studying the LHAM performance for nonlinear PDEs of higher orders, which require a higher homotopy order and finer spatial discretization to obtain accurate solutions.
It was also observed that if the forcing term of any deformation equation has the same decay rate as the homogeneous solution, the forcing term cannot be computed independently in a closed system.
In this case, an extra forcing variable can be added as a linearly independent exponential function multiplied by a monomial time factor, where the order of the monomial scales with the number of successive resonant interactions at the same decay rate.
Furthermore, if the forcing terms of all deformation equations have the same decay rates as the corresponding homogeneous solutions, the final state dimension scales with the orders of the monomial factor.
However, when homogeneous solutions share the same resonant decay rates, the final state dimension is drastically reduced by reusing auxiliary channels with different coupling vectors.
As a result, the dimension required for the secondary linearization scales with the number of distinct resonant rates, which are weighted by their corresponding time-factor monomial orders.

\bibliographystyle{unsrt}
\bibliography{references}

\appendix

\section{\label{sec:Duhamel}Duhamel's Principle for Non-Hermitian Time Evolution}

The nonlinear partial differential equation (PDE) is written with an initial value such as Eq. \eqref{eq:general_PDE}.
In general, non-Hermitian $\mathcal M$ can be decomposed as Hermitian and anti-Hermitian components. Then, the nonunitary time evolution of $\mathcal M$ is defined as $\mathcal G(t) = e^{\mathcal Mt}$, so that $d\mathcal G(t)/dt = \mathcal M\mathcal G(t)$ and $\mathcal G(0)=\mathcal I$.

A new function $\boldsymbol v(\boldsymbol x,t) = \mathcal G(-t)\boldsymbol u(\boldsymbol x,t)$ is introduced such that
\begin{equation}        \label{eq:variation_of_constants}
    \boldsymbol u(\boldsymbol x,t) = \mathcal G(t)\boldsymbol v(\boldsymbol x,t).
\end{equation}
The time derivative of Eq. \eqref{eq:variation_of_constants} is
\begin{equation}
    \label{variation_of_constants_time_derivative}
    \frac{\partial}{\partial t}  \boldsymbol u(\boldsymbol x,t) = \mathcal M \mathcal G(t)\boldsymbol v(\boldsymbol x,t) + \mathcal G(t)\frac{\partial}{\partial t} \boldsymbol v(\boldsymbol x,t).
\end{equation}
Eq. \eqref{variation_of_constants_time_derivative} is substituted by Eq. \eqref{eq:variation_of_constants} such that
\begin{equation}
\label{eq:variation_of_constants_time_derivative_i}
    \frac{\partial}{\partial t} \boldsymbol u(\boldsymbol x,t) = \mathcal M\boldsymbol u(\boldsymbol x,t) + \mathcal G(t)\frac{\partial}{\partial t}\boldsymbol v(\boldsymbol x,t).
\end{equation}
By comparing the right-hand sides of Eq. \eqref{eq:general_PDE} and Eq. \eqref{eq:variation_of_constants_time_derivative_i}, we obtain
\begin{equation}
\label{eq:v_dot}
    \frac{\partial}{\partial t}  \boldsymbol v(\boldsymbol x,t) = \mathcal G(-t) \mathcal N(\boldsymbol u(\boldsymbol x,t)).
\end{equation}
Eq. \eqref{eq:v_dot} is integrated from $0$ to $t$ and the initial condition in Eq. \eqref{eq:general_PDE} is substituted so that
\begin{equation}
\label{eq:v_integral}
    \boldsymbol v(\boldsymbol x,t)= \boldsymbol u_0(\boldsymbol x) +\int_0^t \mathcal G(-\tau) \mathcal N(\boldsymbol u(\boldsymbol x,\tau))d\tau.
\end{equation}
By substituting Eq. \eqref{eq:v_integral} to Eq. \eqref{eq:variation_of_constants}, the exact solution of Eq. \eqref{eq:general_PDE} is obtained in the form of Duhamel's integral such as
\begin{align}
\begin{split}
\label{eq:u_integral}
    \boldsymbol u(\boldsymbol x,t)
    &= \mathcal G(t)\boldsymbol u_0(\boldsymbol x) +\int_0^t \mathcal G(t-\tau) \mathcal N(\boldsymbol u(\boldsymbol x,\tau))d\tau \\
    &= e^{\mathcal Mt} \boldsymbol u_0(\boldsymbol x) +\int_0^t e^{\mathcal M(t-\tau)} \mathcal N(\boldsymbol u(\boldsymbol x,\tau))d\tau.    
\end{split}
\end{align}

\section{\label{sec:proof2}
Proof of Theorem~\ref{thm:stinespring-dilation_error}}

\begin{proof}
Let
\begin{equation}
    \rho(T)
    :=
    e^{\mathcal L T}(\rho_0)
\end{equation}
denote the exact density operator generated by the shifted Lindbladian, and let
\begin{equation}
    \rho_{\tilde n}^{(S)}
    :=
    \left(\mathcal E_{\Delta t}^{(S)}\right)^{\tilde n}(\rho_0),
    \qquad
    \tilde n=\frac{T}{\Delta t},
\end{equation}
denote its order-$S$ Stinespring approximation. Define the accumulated trace-norm error at the $n$th time step by
\begin{equation}
    \varepsilon_n
    :=
    \left\|
        \rho(t_n)-\rho_n^{(S)}
    \right\|_1,
    \qquad
    t_n=n\Delta t.
\end{equation}
The triangle inequality gives
\begin{align}
\varepsilon_n
\leq{}&
\left\|
    e^{\mathcal L\Delta t}
    \left(
        \rho(t_{n-1})-\rho_{n-1}^{(S)}
    \right)
\right\|_1
\nonumber\\
&+
\left\|
    e^{\mathcal L\Delta t}
    \left(\rho_{n-1}^{(S)}\right)
    -
    \mathcal E_{\Delta t}^{(S)}
    \left(\rho_{n-1}^{(S)}\right)
\right\|_1.
\end{align}
Because $e^{\mathcal L\Delta t}$ is a completely positive trace-preserving map, it is contractive on differences of density operators in the trace norm. The local order-$S$ dilation estimate therefore yields
\begin{equation}
    \varepsilon_n
    \leq
    \varepsilon_{n-1}
    +
    C_S
    \mathscr B_{J,\tilde m}^{S+1}
    \Delta t^{S+1},
\end{equation}
where $C_S:=C_{S,1}$ because the LHAM embedding contains one jump operator. Since $\varepsilon_0=0$, iteration over $\tilde n=T/\Delta t$ time steps gives
\begin{equation}
\label{eq:ding_error_sup}
    \left\|
        \rho(T)-\rho_{\tilde n}^{(S)}
    \right\|_1
    \leq
    C_S
    T
    \mathscr B_{J,\tilde m}^{S+1}
    \Delta t^S.
\end{equation}

We next convert this density-operator error into an error in the recovered LHAM solution. Assume that $\boldsymbol y(0)\neq\boldsymbol 0$ and define the normalized initial state by
\begin{equation}
\label{eq:normalized_initial_lham_state}
    \ket{\overline{\boldsymbol y}_0}
    := 
    \frac{1}{\|\boldsymbol y(0)\|_2}
    \sum_{\ell}y_\ell(0)\ket{\ell},
    \qquad
    \left\|
        \ket{\overline{\boldsymbol y}_0}
    \right\|_2=1.
\end{equation}
Define the upper-right-block extraction map by
\begin{equation}
\label{eq:block_extraction_map}
    \mathcal Q_{01}(X)
    :=
    \left(\bra{0}\otimes I\right)
    X
    \left(\ket{1}\otimes I\right).
\end{equation}
In particular,
\begin{equation}
    \rho_{01}(T)
    =
    \mathcal Q_{01}\!\left(\rho(T)\right),
    \qquad
    \rho_{01,\tilde n}^{(S)}
    =
    \mathcal Q_{01}\!\left(\rho_{\tilde n}^{(S)}\right).
\end{equation}

The nondiagonal matrix embedding is constructed such that
\begin{equation}
\label{eq:shifted_block_recovery}
    2\rho_{01}(T)
    \ket{\overline{\boldsymbol y}_0}
    =
    e^{(A-\gamma I)T}
    \ket{\overline{\boldsymbol y}_0}.
\end{equation}
Consequently, the unshifted lifted LHAM vector is recovered as
\begin{equation}
\label{eq:exact_lifted_recovery}
    \boldsymbol y(T)
    =
    2\|\boldsymbol y(0)\|_2
    e^{\gamma T}
    \rho_{01}(T)
    \ket{\overline{\boldsymbol y}_0}.
\end{equation}
Analogously, define the approximate lifted vector by
\begin{equation}
\label{eq:approximate_lifted_recovery}
    \widehat{\boldsymbol y}_{S,\Delta t}(T)
    :=
    2\|\boldsymbol y(0)\|_2
    e^{\gamma T}
    \rho_{01,\tilde n}^{(S)}
    \ket{\overline{\boldsymbol y}_0}.
\end{equation}

Let $\mathcal R_{J,\tilde m}$ denote the reconstruction map from the lifted LHAM vector to the truncated homotopy-series solution:
\begin{equation}
    \boldsymbol u_{J,\tilde m}(T)
    =
    \mathcal R_{J,\tilde m}\boldsymbol y(T),
    \qquad
    \widehat{\boldsymbol u}_{J,\tilde m,S,\Delta t}(T)
    =
    \mathcal R_{J,\tilde m}
    \widehat{\boldsymbol y}_{S,\Delta t}(T).
\end{equation}
It follows that
\begin{align}
&
\left\|
    \boldsymbol u_{J,\tilde m}(T)
    -
    \widehat{\boldsymbol u}_{J,\tilde m,S,\Delta t}(T)
\right\|_2
\nonumber\\
&\leq
2\|\boldsymbol y(0)\|_2
e^{\gamma T}
\|\mathcal R_{J,\tilde m}\|_2
\left\|
    \mathcal Q_{01}
    \left(
        \rho(T)-\rho_{\tilde n}^{(S)}
    \right)
    \ket{\overline{\boldsymbol y}_0}
\right\|_2
\nonumber\\
&\leq
2\|\boldsymbol y(0)\|_2
e^{\gamma T}
\|\mathcal R_{J,\tilde m}\|_2
\left\|
    \mathcal Q_{01}
    \left(
        \rho(T)-\rho_{\tilde n}^{(S)}
    \right)
\right\|_2.
\label{eq:dilation_recovery_intermediate}
\end{align}
For any matrix $X$, the block-extraction map satisfies
\begin{align}
    \|\mathcal Q_{01}(X)\|_2
    &\leq
    \|\mathcal Q_{01}(X)\|_1
    \nonumber\\
    &=
    \left\|
        (\bra{0}\otimes I)
        X
        (\ket{1}\otimes I)
    \right\|_1
    \nonumber\\
    &\leq
    \|\bra{0}\otimes I\|_2
    \|X\|_1
    \|\ket{1}\otimes I\|_2
    \nonumber\\
    &=
    \|X\|_1.
\label{eq:block_extraction_contraction}
\end{align}
Therefore,
\begin{equation}
\label{eq:recovery_trace_norm_bound}
\begin{aligned}
&
\left\|
    \boldsymbol u_{J,\tilde m}(T)
    -
    \widehat{\boldsymbol u}_{J,\tilde m,S,\Delta t}(T)
\right\|_2
\\
&\qquad\leq
2\|\boldsymbol y(0)\|_2
e^{\gamma T}
\|\mathcal R_{J,\tilde m}\|_2
\left\|
    \rho(T)-\rho_{\tilde n}^{(S)}
\right\|_1.
\end{aligned}
\end{equation}

For the LHAM reconstruction,
\begin{equation}
    \mathcal R_{J,\tilde m}
    =
    [\,I\;\;0\;\;I\;\;0\;\cdots\;0\;\;I\,],
\end{equation}
where there are exactly $\tilde m+1$ identity blocks. Hence,
\begin{equation}
    \mathcal R_{J,\tilde m}
    \mathcal R_{J,\tilde m}^{\dagger}
    =
    (\tilde m+1)I,
\end{equation}
and thus
\begin{equation}
\label{eq:reconstruction_operator_norm}
    \|\mathcal R_{J,\tilde m}\|_2
    =
    \sqrt{\tilde m+1}.
\end{equation}
Substituting Eqs.~\eqref{eq:ding_error_sup} and \eqref{eq:reconstruction_operator_norm} into Eq.~\eqref{eq:recovery_trace_norm_bound} yields
\begin{equation}
\begin{aligned}
&
\left\|
    \boldsymbol u_{J,\tilde m}(T)
    -
    \widehat{\boldsymbol u}_{J,\tilde m,S,\Delta t}(T)
\right\|_2
\\
&\qquad\leq
2C_S\|\boldsymbol y(0)\|_2
\sqrt{\tilde m+1}\,
e^{\gamma T}
T\mathscr B_{J,\tilde m}^{S+1}
\Delta t^S,
\end{aligned}
\end{equation}
which proves Eq.~\eqref{eq:dilation_error}.
\end{proof}

\section{\label{sec:proof3}Proof of Theorem \ref{thm:homotopy_dilation_balance}}

\begin{proof}
Insert the Fourier-truncated solution $\boldsymbol u_J(T)$ and the order-$\tilde m$ homotopy approximation $\boldsymbol u_{J,\tilde m}(T)$. 
The triangle inequality gives
\begin{equation}
\label{eq:total_error_triangle_decomposition}
\begin{aligned}
    \|\boldsymbol u(T)-\widehat{\boldsymbol{u}}_{J,\tilde{m},S,\Delta t}(T)\|_{L^2(\Omega)}
    \le\;&
    \left\|\boldsymbol u(T)-\boldsymbol u_J(T)\right\|_{L^2(\Omega)}
    \\
    &+
    \|\boldsymbol u_J(T)-\boldsymbol u_{J,\tilde{m}}(T)\|_{L^2(\Omega)}
    \\
    &+
    \|\boldsymbol u_{J,\tilde{m}}(T)-\widehat{\boldsymbol{u}}_{J,\tilde{m},S,\Delta t}(T)\|_{L^2(\Omega)}.
\end{aligned}
\end{equation}

The first term is the Fourier truncation error. Let $\Omega=\mathbb T^{d_x}$ be a $d_x$-dimensional periodic domain, where $\mathbb T:=\mathbb R/(2\pi\mathbb Z)$ and $d_x$ is the number of spatial dimensions. Equivalently, $\Omega$ may be identified with $[0,2\pi]^{d_x}$ subject to periodic boundary conditions. Assume that $
    \boldsymbol u
    \in
    \mathscr C\!\left([0,T];H^f(\Omega)\right)$,
where $f\geq 0$ and $\mathscr C([0,T];H^f(\Omega))$ denotes the space of continuous mappings from $[0,T]$ into the Sobolev space $H^f(\Omega)$. Thus, $\boldsymbol u(t)$ has spatial Sobolev regularity of order $f$ and depends continuously on $t$ in the $H^f$-norm. 
Since the mapping $t\mapsto\|\boldsymbol u(t)\|_{H^f(\Omega)}$ is continuous on the compact interval $[0,T]$, it attains a finite maximum. Therefore,
\begin{equation}
    \|\boldsymbol u\|_{L^\infty(0,T;H^f(\Omega))}
    =
    \max_{0\leq t\leq T}
    \|\boldsymbol u(t)\|_{H^f(\Omega)}
    <\infty.
\end{equation}

Let $\boldsymbol j=(j_1,\ldots,j_{d_x})\in\mathbb Z^{d_x}$ denote a Fourier multi-index, and, for $J\in\mathbb N$, define the tensor-product Fourier space $
    \mathbb B_J
    :=
    \operatorname{span}
    \left\{
        e^{i\boldsymbol j\cdot\boldsymbol x}:
        |j_k|\leq J,\quad k=1,\ldots,d_x
    \right\}$.
Let $P_J:L^2(\Omega)\rightarrow\mathbb B_J$ denote the $L^2$-orthogonal projection, defined by
\begin{equation}
    \left\langle
        b-P_Jb,a_J
    \right\rangle_{L^2(\Omega)}
    =0
    \qquad
    \forall a_J\in\mathbb B_J.
\end{equation}
For vector-valued functions, $P_J$ is applied componentwise.
The standard Fourier projection estimate gives
\cite{shen2011spectral,Hesthaven_Gottlieb_Gottlieb_2007}
\begin{equation}
\begin{aligned}
    \left\|
        \boldsymbol u(T)-P_J\boldsymbol u(T)
    \right\|_{L^2(\Omega)}
    &\leq
    C J^{-f}
    \left\|\boldsymbol u(T)\right\|_{H^f(\Omega)}
    \\
    &\leq
    C J^{-f}
    \left\|\boldsymbol u\right\|_
    {L^\infty(0,T;H^f(\Omega))},
\end{aligned}
\end{equation}
where $C>0$ is independent of $J$. Defining $
    \boldsymbol u_J(t):=P_J\boldsymbol u(t)$
and $
    C_f
    :=
    C
    \left\|\boldsymbol u\right\|_
    {L^\infty(0,T;H^f(\Omega))}$,
we obtain
\begin{equation}
\label{eq:fourier_error_in_total_proof}
    \left\|
        \boldsymbol u(T)-\boldsymbol u_J(T)
    \right\|_{L^2(\Omega)}
    \leq
    C_fJ^{-f}.
\end{equation}

The second term is the homotopy truncation error.
Assume that the Fourier-truncated solution admits the convergent homotopy-series solution
\begin{equation}
    \boldsymbol u_J(T)
    =
    \sum_{m=0}^{\infty}\boldsymbol u_J^{(m)}(T),
\end{equation}
and define its order-$\tilde m$ truncation by
\begin{equation}
    \boldsymbol u_{J,\tilde m}(T)
    :=
    \sum_{m=0}^{\tilde m}\boldsymbol u_J^{(m)}(T).
\end{equation}
Suppose that the homotopy coefficients satisfy the geometric-decay condition
\begin{equation}
\label{eq:homotopy_geometric_decay}
    \left\|
        \boldsymbol u_J^{(m+1)}(T)
    \right\|_{L^2(\Omega)}
    \leq
    \theta_J
    \left\|
        \boldsymbol u_J^{(m)}(T)
    \right\|_{L^2(\Omega)},
    \qquad
    m\geq0,
\end{equation}
for some contraction factor $\theta_J\in(0,1)$. Then the triangle inequality and the geometric-series formula give
\begin{align}
\left\|
    \boldsymbol u_J(T)-\boldsymbol u_{J,\tilde m}(T)
\right\|_{L^2(\Omega)}
\nonumber
&\leq
\sum_{m=\tilde m+1}^{\infty}
\left\|
    \boldsymbol u_J^{(m)}(T)
\right\|_{L^2(\Omega)}
\nonumber\\
&\leq
\left\|
    \boldsymbol u_J^{(0)}(T)
\right\|_{L^2(\Omega)}
\sum_{m=\tilde m+1}^{\infty}\theta_J^m
\nonumber\\
&=
\left\|
    \boldsymbol u_J^{(0)}(T)
\right\|_{L^2(\Omega)}
\frac{\theta_J^{\tilde m+1}}{1-\theta_J}.
\label{eq:homotopy_error_in_total_proof}
\end{align}

The third term is Stinespring dilation error. Theorem~\ref{thm:stinespring-dilation_error} gives
\begin{align}
&
\left\|
    \boldsymbol u_{J,\tilde m}(T)
    -
    \widehat{\boldsymbol u}_{J,\tilde m,S,\Delta t}(T)
\right\|_{L^2(\Omega)}
\nonumber\\
&\leq
2C_S\|\boldsymbol y(0)\|_2
\sqrt{\tilde m+1}\,
e^{\gamma T}
T\mathscr B_{J,\tilde m}^{S+1}
\Delta t^S.
\label{eq:dilation_error_in_total_proof}
\end{align}
Here the Euclidean norm of the Fourier-coefficient error is identified with the spatial $L^2(\Omega)$-norm through Parseval's identity, using an orthonormal Fourier basis.

Combining
Eqs.~\eqref{eq:fourier_error_in_total_proof}, \eqref{eq:homotopy_error_in_total_proof}, and \eqref{eq:dilation_error_in_total_proof} with Eq.~\eqref{eq:total_error_triangle_decomposition} gives
\begin{align}
\left\|
    \boldsymbol u(T)
    -
    \widehat{\boldsymbol u}_{J,\tilde m,S,\Delta t}(T)
\right\|_{L^2(\Omega)}
\nonumber
\leq\;&
C_fJ^{-f} \\
&+
\left\|
    \boldsymbol u_J^{(0)}(T)
\right\|_{L^2(\Omega)}
\frac{\theta_J^{\tilde m+1}}{1-\theta_J}
\nonumber\\
&+
2C_S\|\boldsymbol y(0)\|_2
\sqrt{\tilde m+1}\,
e^{\gamma T}
T\mathscr B_{J,\tilde m}^{S+1}
\Delta t^S,
\end{align}
which proves Eq.~\eqref{eq:lham_homotopy_dilation_error}.
\end{proof}

\section{\label{sec:Burgers}LHAM Derivation of Burgers' Equation}

The orthonormal Fourier basis in the periodic spatial domain $x \in [0,2\pi)$ is chosen as
\begin{equation}
    \varphi_j(x)=
    e^{i j x},
    \qquad j\in\mathbb Z,
\end{equation}
and truncated to $|j|\le J$
such that $\boldsymbol{c}(t)=(c_{-J}(t),\ldots,c_J(t)) \in\mathbb C^{2J+1}$.

By projecting the terms in Eq. \eqref{eq:general_PDE} onto $\mathrm{span}\{\varphi_j\}_{j=-J}^J$, an ordinary differential equation (ODE) system is obtained as Eq. \eqref{eq:ode_coeff}.
Thus, Eq. \eqref{eq:ode_coeff} is defined with $M =[M_{kj}]\in\mathbb{R}^{(2J+1) \times (2J+1)} $ and $\boldsymbol{N}(\boldsymbol{c}(t)) = (N_{-J}(\boldsymbol{c}(t)), \ldots, N_{J}(\boldsymbol{c}(t)))$.

The linear part is diagonalized as
\begin{equation}
    \ell_{kj}
    = \braket{\varphi_k, \mathcal{M}\varphi_j} 
    =\nu\langle \varphi_k,\frac{\partial^2}{\partial x^2}\varphi_j\rangle
    =-\nu j^2\,\delta_{kj}.
    \label{eq:burgers_A_diag}
\end{equation}
The nonlinear term in Eq. \eqref{eq:burgers_pde} is 
projected and truncated to $|j|\le J$, then the $j$th-component of $\boldsymbol{N}(\boldsymbol{c}(t))$ is defined as
\begin{equation}
    N_j(\boldsymbol{c})
    =
    -\sum_{\substack{a+b=j\\ |a|,|b|\le J}}
    (ib)\,c_a\,c_b,
    \qquad |j|\le J.
    \label{eq:burgers_F_conv}
\end{equation}

By introducing homotopy-Maclaurin series in Eq. \eqref{eq:homotopy_Maclaurin_series_f} with the initial conditions in Eq. \eqref{eq:init_conditions_f}, the $m$th-order deformation equation is obtained as Eq. \eqref{eq:forced_linear_cm} for $m\ge1$.

The linear time evolution $e^{M_{jj} t}=e^{-\nu j^2 t}$ corresponds to mode-wise exponential damping, as expected for viscosity. The Lindbladian homotopy analysis method (LHAM) workflow is then classically construct matrix $A$ and evaluate Eq. \eqref{eq:augmented_linear_system_f} using the Lindbladian.

For the initial condition $u(x,0)=\sin x$, $\sin x = (e^{ix}-e^{-ix})/2i$ is expanded with the Fourier basis, where coefficients are $c_{-1}^{(0)}(0) = {i}/{2}$, $c_{1}^{(0)}(0) = -i/2$, and zeros for all others.
From Eq. \eqref{eq:c0_hom}, the time evolution for $\boldsymbol{c}^{(0)}(t)$ is therefore
\begin{equation}
c_j^{(0)}(t)=e^{-\nu j^2 t}c_j^{(0)}(0).
\end{equation}
From Eq. \eqref{eq:Galerkin}, the zeroth-order solution is 
\begin{equation}
\label{eq:zero_solution}
\tilde{u}^{(0)}(x,t)=e^{-\nu t}\sin x.
\end{equation}
The forcing term in the first-order deformation equation is
\begin{equation}
\label{eq:1th_sorcing}
f^{(0)} = -\tilde u^{(0)} \frac{\partial}{\partial x} \tilde u^{(0)} .
\end{equation}
Substituting the zeroth-order solution Eq. \eqref{eq:zero_solution} in Eq. \eqref{eq:1th_sorcing} yields
\begin{equation}
f^{(0)}(x,t)=-\frac{1}{2}e^{-2\nu t}\sin(2x).
\end{equation}
In the reciprocal space, this produces only $j=\pm2$ modes
\begin{align}
f_{-2}^{(0)}(t) = -\frac{i}{4}e^{-2\nu t}, \qquad
f_{2}^{(0)}(t) = \frac{i}{4}e^{-2\nu t}.
\end{align}

Thus, the forcing vector can be written as Eq. \eqref{eq:single_channel_force_f}.
To remove explicit time dependence we introduce the auxiliary variable
$z^{(1)}(t)= e^{-2\nu t}$ satisfying $d{z}^{(1)}(t)/dt = -2\nu z^{(1)}(t)$.
The lifted state defined as Eq. \eqref{eq:mu_def_f} gives the autonomous linear system with block matrix
\begin{equation}
A =
\begin{bmatrix}
M & 0 & 0 & \cdots \\
0 & -\lambda^{(1)} & 0 & \cdots \\
0 & \boldsymbol{v}^{(1)}  & M & \cdots \\
\vdots & \vdots & \vdots & \ddots
\end{bmatrix}.
\end{equation}
where $\lambda^{(1)} = 2\nu$, $v^{(1)}_{-2} = -i/4$, $v^{(1)}_{2} = i/4$, and $v^{(1)}_{j} = 0$ for the rest of the $j$-basis.

\section{\label{sec:MHD}LHAM Derivation of Magnetohydrodynamics}

A vector field $\boldsymbol Y(\boldsymbol x,t) = (\omega(\boldsymbol x,t), \xi(\boldsymbol x,t))$ is defined so that the reduced magnetohydrodynamics (MHD) system Eq. \eqref{eq:rmhd_pde} is written as
\begin{equation}
    \frac{\partial}{\partial t} \boldsymbol Y(\boldsymbol{x},t)
    =\mathcal{M}\boldsymbol Y(\boldsymbol{x},t)
    +\mathcal{N}(\boldsymbol Y(\boldsymbol{x},t),\boldsymbol Y(\boldsymbol{x},t)),
\end{equation}
where the linear operator is
\begin{equation}
    \mathcal{M}\begin{bmatrix}
    \omega \\ \xi\end{bmatrix} 
    = \begin{bmatrix}\nu \nabla^2\omega \\ \eta\nabla^2\xi\end{bmatrix},
\end{equation}
and the bilinear nonlinear operator is defined as
\begin{equation}
\label{eq:MHD_nonlinear_op}
\begin{split}
&\mathcal{N}(\boldsymbol Y^{(a)},\boldsymbol Y^{(b)}) \\
&=
\begin{bmatrix}
-\left(
\frac{\partial \phi^{(a)}}{\partial x}\frac{\partial \omega^{(b)}}{\partial y}
-
\frac{\partial \phi^{(a)}}{\partial y}\frac{\partial \omega^{(b)}}{\partial x}
\right)
+
\left(
\frac{\partial \xi^{(a)}}{\partial x}\frac{\partial \zeta^{(b)}}{\partial y}
-
\frac{\partial \xi^{(a)}}{\partial y}\frac{\partial \zeta^{(b)}}{\partial x}
\right)
\\[1em]
-\left(
\frac{\partial \phi^{(a)}}{\partial x}\frac{\partial \xi^{(b)}}{\partial y}
-
\frac{\partial \phi^{(a)}}{\partial y}\frac{\partial \xi^{(b)}}{\partial x}
\right)
\end{bmatrix}.
\end{split}
\end{equation}
with $\phi^{(a)}=-\nabla^{-2}\omega^{(a)}$ and $\zeta^{(b)}=-\nabla^2\xi^{(b)}$.

The orthonormal Fourier basis on the periodic spatial domain $\boldsymbol x=(x,y)\in[0,2\pi)^2$ is chosen as
\begin{equation}
    \varphi_{\boldsymbol j}(\boldsymbol x)
    =
    e^{i\boldsymbol j\cdot \boldsymbol x},
    \qquad
    \boldsymbol j=(j_x,j_y)\in\mathbb Z^2,
\end{equation}
and truncated to $|j_x|,|j_y|\le J$, such that
\begin{equation}
    \omega(\boldsymbol x,t)
    \approx
    \tilde\omega(\boldsymbol x,t)
    :=
    \sum_{|j_x|,|j_y|\le J}
    \omega_{\boldsymbol j}(t)\,
    \varphi_{\boldsymbol j}(\boldsymbol x),
\end{equation}
and
\begin{equation}
    \xi(\boldsymbol x,t)
    \approx
    \tilde\xi(\boldsymbol x,t)
    :=
    \sum_{|j_x|,|j_y|\le J}
    \xi_{\boldsymbol j}(t)\,
    \varphi_{\boldsymbol j}(\boldsymbol x).
\end{equation}
Thus, the coefficient vector is defined as
\begin{equation}
    \boldsymbol c(t)
    =
    \bigl(\boldsymbol\omega(t),\boldsymbol\xi(t)\bigr)
    \in
    \mathbb C^{2(2J+1)^2},
\end{equation}
where $\boldsymbol\omega(t)$ and $\boldsymbol\xi(t)$ are from the Fourier coefficients $\omega_{\boldsymbol j}(t)$ and $\xi_{\boldsymbol j}(t)$, respectively.

By projecting the terms in Eq. \eqref{eq:general_PDE} onto $\mathrm{span}\{\varphi_{\boldsymbol j}\}_{|j_x|,|j_y|\le J}$, an ODE system is obtained as Eq. \eqref{eq:ode_coeff}.
Thus, Eq. \eqref{eq:ode_coeff} is defined with the linear matrix $M$ associated with $\mathcal M$ and the nonlinear vector $\boldsymbol N(\boldsymbol c(t))$.

The linear part is diagonalized in Fourier space. Since
\begin{equation}
    \nabla^2 \varphi_{\boldsymbol j}
    =
    -|\boldsymbol j|^2 \varphi_{\boldsymbol j},
    \qquad
    |\boldsymbol j|^2=j_x^2+j_y^2,
\end{equation}
the linear evolution of each Fourier mode is
\begin{equation}
    \partial_t \omega_{\boldsymbol j}
    =
    -\nu |\boldsymbol j|^2 \omega_{\boldsymbol j},
    \qquad
    \partial_t \xi_{\boldsymbol j}
    =
    -\eta |\boldsymbol j|^2 \xi_{\boldsymbol j}.
\end{equation}
Therefore, the matrix $M$ is block diagonal,
\begin{equation}
    M
    =
    \begin{bmatrix}
    M_\omega & 0\\
    0 & M_\xi
    \end{bmatrix},
\end{equation}
where
\begin{equation}
    (M_\omega)_{\boldsymbol k\boldsymbol j}
    =
    -\nu |\boldsymbol j|^2 \delta_{\boldsymbol k\boldsymbol j},
    \qquad
    (M_\xi)_{\boldsymbol k\boldsymbol j}
    =
    -\eta |\boldsymbol j|^2 \delta_{\boldsymbol k\boldsymbol j}.
\end{equation}

The nonlinear terms in Eq. \eqref{eq:rmhd_pde} are projected and truncated to $|j_x|,|j_y|\le J$.
In Fourier space, each Poisson bracket becomes a modal convolution. For two fields
$f(\boldsymbol x)=\sum_{\boldsymbol a} f_{\boldsymbol a} e^{i\boldsymbol a\cdot\boldsymbol x}$
and
$g(\boldsymbol x)=\sum_{\boldsymbol b} g_{\boldsymbol b} e^{i\boldsymbol b\cdot\boldsymbol x}$,
the $\boldsymbol j$th Fourier coefficient of $(\frac{\partial f}{\partial x}\frac{\partial g}{\partial y} - \frac{\partial f}{\partial y}\frac{\partial g}{\partial x})$ is
\begin{equation}
    \bigl(\frac{\partial f}{\partial x}\frac{\partial g}{\partial y} - \frac{\partial f}{\partial y}\frac{\partial g}{\partial x}\bigr)_{\boldsymbol j}
    =
    -\sum_{\boldsymbol a+\boldsymbol b=\boldsymbol j}
    (a_x b_y-a_y b_x)\,
    f_{\boldsymbol a}\,g_{\boldsymbol b}.
    \label{eq:mhd_bracket_conv}
\end{equation}
Using
\begin{equation}
    \phi_{\boldsymbol a}
    =
    \frac{\omega_{\boldsymbol a}}{|\boldsymbol a|^2},
    \qquad
    \zeta_{\boldsymbol b}
    =
    |\boldsymbol b|^2 \xi_{\boldsymbol b},
    \qquad
    (\boldsymbol a\neq \boldsymbol 0,\ \boldsymbol b\neq \boldsymbol 0),
\end{equation}
the $\boldsymbol j$th nonlinear coefficient in the vorticity equation is
\begin{equation}
\begin{split}
    & N^{(\omega)}_{\boldsymbol j}(\boldsymbol c) \\
    &=
    -\sum_{\boldsymbol a+\boldsymbol b=\boldsymbol j}
    (a_x b_y-a_y b_x)\,
    \phi_{\boldsymbol a}\,\omega_{\boldsymbol b}
    +
    \sum_{\boldsymbol a+\boldsymbol b=\boldsymbol j}
    (a_x b_y-a_y b_x)\,
    \xi_{\boldsymbol a}\,\zeta_{\boldsymbol b},
\end{split}
\end{equation}
and the $\boldsymbol j$th nonlinear coefficient in the magnetic potential equation is
\begin{equation}
    N^{(\xi)}_{\boldsymbol j}(\boldsymbol c)
    =
    -\sum_{\boldsymbol a+\boldsymbol b=\boldsymbol j}
    (a_x b_y-a_y b_x)\,
    \phi_{\boldsymbol a}\,\xi_{\boldsymbol b}.
\end{equation}
Thus, the nonlinear vector is written as
\begin{equation}
    \boldsymbol N(\boldsymbol c)
    =
    \bigl(
    N^{(\omega)}_{\boldsymbol j}(\boldsymbol c),
    N^{(\xi)}_{\boldsymbol j}(\boldsymbol c)
    \bigr)_{|j_x|,|j_y|\le J}.
\end{equation}

By introducing the homotopy-Maclaurin series in Eq. \eqref{eq:homotopy_Maclaurin_series_f} with the initial conditions in Eq. \eqref{eq:init_conditions_f}, the $m$th-order deformation equation is obtained as Eq. \eqref{eq:forced_linear_cm} for $m\ge 1$,
where the forcing vector is
\begin{equation}
\label{eq:MHD_0th_forcing}
    \boldsymbol f^{(m-1)}(t)
    =
    \sum_{\alpha=0}^{m-1}
    \boldsymbol N\!\left(
    \boldsymbol c^{(\alpha)}(t),
    \boldsymbol c^{(m-1-\alpha)}(t)
    \right).
\end{equation}
Equivalently,
\begin{equation}
    \boldsymbol N\bigl(\boldsymbol\Phi(t;q),\boldsymbol\Phi(t;q)\bigr)
    =
    \sum_{m\ge 0}\boldsymbol S^{(m)}(t)\,q^m,
\end{equation}
where $\boldsymbol S^{(m)}(t)$ is the coefficient of $q^m$ obtained from the bilinear convolution sums.

The linear time evolution is defined as $e^{Mt}$.
Since $M$ is time-independent and diagonal in Fourier space, the propagator acts mode-wise as
\begin{equation}
\label{eq:MHD_init_linear_diff_eq}
    \omega^{(0)}_{\boldsymbol j}(t)
    =
    e^{-\nu |\boldsymbol j|^2 t}\,\omega^{(0)}_{\boldsymbol j}(0),
    \qquad
    \xi^{(0)}_{\boldsymbol j}(t)
    =
    e^{-\eta |\boldsymbol j|^2 t}\,\xi^{(0)}_{\boldsymbol j}(0),
\end{equation}
where $|\boldsymbol j|^2=j_x^2+j_y^2$.
Thus, the linear time evolution corresponds to mode-wise exponential damping of the vorticity and magnetic potential modes, as expected from viscosity and resistivity.

For the initial condition Eq. \eqref{eq:MHD_IC}, the Fourier basis expansion is obtained from
\begin{equation}
    \varphi_{\boldsymbol j}(x,y)=e^{i(j_x x+j_y y)},
    \qquad
    \boldsymbol j=(j_x,j_y)\in\mathbb Z^2.
\end{equation}
Using
\begin{equation}
    \sin x=\frac{e^{ix}-e^{-ix}}{2i},
    \qquad
    \sin(x-y)=\frac{e^{i(x-y)}-e^{-i(x-y)}}{2i},
\end{equation}
and
\begin{equation}
    \cos y=\frac{e^{iy}+e^{-iy}}{2},
    \qquad
    \cos(x+y)=\frac{e^{i(x+y)}+e^{-i(x+y)}}{2},
\end{equation}
the nonzero Fourier coefficients of the initial vorticity are
\begin{align}
    \omega^{(0)}_{(-1,0)}(0)&=\frac{i}{2}, &
    \omega^{(0)}_{(1,0)}(0)&=-\frac{i}{2}, \\
    \omega^{(0)}_{(-1,1)}(0)&=\frac{i}{4}, &
    \omega^{(0)}_{(1,-1)}(0)&=-\frac{i}{4},
\end{align}
and zeros for all other modes.
Similarly, the nonzero Fourier coefficients of the initial magnetic potential are
\begin{align}
    \xi^{(0)}_{(0,-1)}(0)&=\frac12, &
    \xi^{(0)}_{(0,1)}(0)&=\frac12, \\
    \xi^{(0)}_{(-1,-1)}(0)&=\frac18, &
    \xi^{(0)}_{(1,1)}(0)&=\frac18,
\end{align}
and zeros for all other modes.

From the initial linear differential equation,
\begin{equation}
    \frac{\partial}{\partial t} \omega^{(0)}=\nu \nabla^2 \omega^{(0)},
    \qquad
    \frac{\partial}{\partial t} \xi^{(0)}=\eta \nabla^2 \xi^{(0)},
\end{equation}
each Fourier mode evolves independently as Eq. \eqref{eq:MHD_init_linear_diff_eq}

Therefore, the zeroth-order solutions in physical space are
\begin{equation}
\label{eq:mhd_zero_omega}
    \tilde\omega^{(0)}(x,y,t)
    =
    e^{-\nu t}\sin x
    +
    \frac12 e^{-2\nu t}\sin(x-y),
\end{equation}
and
\begin{equation}
\label{eq:mhd_zero_xi}
    \tilde\xi^{(0)}(x,y,t)
    =
    e^{-\eta t}\cos y
    +
    \frac14 e^{-2\eta t}\cos(x+y).
\end{equation}
Here the modes $(\pm1,0)$ and $(0,\pm1)$ decay at rates $\nu$ and $\eta$, respectively, while the diagonal modes $(\pm1,\mp1)$ and $(\pm1,\pm1)$ decay at rates $2\nu$ and $2\eta$, respectively.

Using $\phi^{(0)}=-\nabla^{-2}\omega^{(0)}$ and $\zeta^{(0)}=-\nabla^2\xi^{(0)}$, one obtains
\begin{equation}
    \tilde \phi^{(0)}(x,y,t)
    =
    e^{-\nu t}\sin x
    +
    \frac14 e^{-2\nu t}\sin(x-y),
\end{equation}
and
\begin{equation}
    \tilde \zeta^{(0)}(x,y,t)
    =
    e^{-\eta t}\cos y
    +
    \frac12 e^{-2\eta t}\cos(x+y).
\end{equation}

The first-order forcing term is then obtained from the bilinear nonlinear operator in Eq. \eqref{eq:MHD_0th_forcing} evaluated at the zeroth-order solution.
Thus, the first-order deformation equation is a forced linear system whose source term is entirely determined by the zeroth-order modes.

Since each zeroth-order Fourier mode evolves exponentially,
\begin{equation}
    \omega_{\boldsymbol j}^{(0)}(t)
    \sim e^{-\nu |\boldsymbol j|^2 t},
    \qquad
    \xi_{\boldsymbol j}^{(0)}(t)
    \sim e^{-\eta |\boldsymbol j|^2 t},
\end{equation}
every product appearing in the nonlinear operator produces another exponential factor. 
Thus, the forcing vector at order $m$ can be written as Eq. \eqref{eq:single_channel_force_f} with auxiliary variables in Eq. \eqref{eq:aux_variables}, and $\lambda^{(m)}_{r}$ are obtained by summing the decay rates of the lower-order modes appearing in the corresponding nonlinear product. For example, if a forcing contribution contains the product $e^{-\lambda_a t}e^{-\lambda_b t}$, then $\lambda^{(m)}_{r}=\lambda_a+\lambda_b$.

Once the forcing vector $\boldsymbol f^{(m-1)}(t)$ is evaluated, each component of Eq. \eqref{eq:forced_linear_cm} evolves independently.
For a single Fourier mode $\boldsymbol j$, let $c(t)$ denote either $\omega_{\boldsymbol j}^{(m)}(t)$ or $\xi_{\boldsymbol j}^{(m)}(t)$.
For a single Fourier mode and a single forcing channel, Eq. \eqref{eq:forced_linear_cm} reduces component-wise to
\begin{equation}
\label{eq:MHD_componentwise}
    \frac{d}{dt}c(t)
    =
    -\sigma c(t)+v e^{-\lambda t}, \qquad
    c(0)=0,
\end{equation}
where $\sigma\in\{\nu |\boldsymbol k|^2,\eta |\boldsymbol k|^2\}$ and $\lambda$ is one of the forcing decay rates $\lambda^{(m)}_r$.
If $\sigma\neq \lambda$, the solution is
\begin{equation}
    c(t)
    =
    \frac{v}{\sigma-\lambda}
    \left(
    e^{-\lambda t}-e^{-\sigma t}
    \right).
\end{equation}

To remove explicit time dependence, the forcing terms are represented through auxiliary variables satisfying Eq. \eqref{eq:aux_variables}.
If the forcing at a given order is written as Eq. \eqref{eq:single_channel_force_f} with Eq. \eqref{eq:zm_evolution},
then the lifted state defined as Eq. \eqref{eq:mu_def_f} with $\boldsymbol{c}(t)$ gives an linear homogeneous autonomous system with lower block triangular matrix Eq. \eqref{eq:A_block_form}, or explicitly for one auxilliary variable with Eq. \eqref{eq:MHD_componentwise}
\begin{equation}
A =
\begin{bmatrix}
\ddots & \vdots & \vdots & \vdots & \iddots \\
\cdots & M & 0 & 0 & \cdots \\
\cdots & 0 & -\lambda & 0 & \cdots \\
\cdots & 0 & v  & M & \cdots \\
\iddots & \vdots & \vdots & \vdots & \ddots
\end{bmatrix},
\end{equation}
where $M=-\sigma$ and each block $V^{(m)}$ couples the auxiliary forcing variables to the $m$th-order modal coefficients.

Therefore, similarly to the Burgers' equation, the LHAM workflow for reduced MHD is to classically construct the block matrix $A$ from the Fourier basis-projected linear operator and homotopy forcing vectors, and then evaluate Eq. \eqref{eq:augmented_linear_system_f} using the Lindblad dynamics.

\nocite{*}



\end{document}